\documentclass{my_aims}

\usepackage{amsmath}
\usepackage{paralist}
\usepackage{changepage}
\usepackage{subfig}
\usepackage{graphics} 
\usepackage{epsfig} 
\usepackage{graphicx}  

\usepackage{url}

\def\ppages{X--XX}

\theoremstyle{definition}

\title[Ebola Model and Optimal Control]{Ebola Model and Optimal Control\\ 
with Vaccination Constraints}

\author[I. Area, F. Nda\"{\i}rou, J. J. Nieto, C. J. Silva and D. F. M. Torres]{}

\subjclass{Primary: 49J15, 92D30; Secondary: 34C60, 49N90.}

\keywords{Ebola virus, mathematical modelling, transmission of Ebola,
control of the spread of the Ebola disease,
optimal control with vaccination constraints,
vaccination scenarios.}

\email{area@uvigo.es}
\email{faical.ndairou@aims-cameroon.org}
\email{juanjose.nieto.roig@usc.es}
\email{cjoaosilva@ua.pt}
\email{delfim@ua.pt}

\thanks{$^*$Corresponding author: delfim@ua.pt (D. F. M. Torres)}

\begin{document}

\maketitle

\centerline{\scshape Iv\'an Area}
\medskip
{\footnotesize \centerline{Departamento de Matem\'{a}tica Aplicada II, 
E. E. Aeron\'autica e do Espazo}
\centerline{Campus As Lagoas, Universidade de Vigo, 32004 Ourense, Spain}
} 


\medskip

\centerline{\scshape Fa\"{\i}\c{c}al Nda\"{\i}rou}
\medskip
{\footnotesize
\centerline{African Institute for Mathematical Sciences (AIMS--Cameroon)}
\centerline{P.O. Box 608, Limbe Crystal Gardens, South West Region, Cameroon}
}


\medskip

\centerline{\scshape Juan J. Nieto}
\medskip
{\footnotesize
\centerline{Departamento de An\'alise Matem\'atica, Estat\'{\i}stica e Optimizaci\'on}
\centerline{Facultade de Matem\'{a}ticas, Universidade de Santiago de Compostela}
\centerline{15782 Santiago de Compostela, Spain}
}


\medskip

\centerline{\scshape Cristiana J. Silva and Delfim F. M. Torres$^*$}
\medskip
{\footnotesize
\centerline{Center for Research and Development in Mathematics and Applications (CIDMA)}
\centerline{Department of Mathematics, University of Aveiro, 3810-193 Aveiro, Portugal}
} 


\begin{abstract}
The Ebola virus disease is a severe viral haemorrhagic fever 
syndrome caused by Ebola virus. This disease is transmitted by direct contact 
with the body fluids of an infected person and objects contaminated with virus 
or infected animals, with a death rate close to 90\% in humans. Recently, 
some mathematical models have been presented to analyse the spread of the 2014 
Ebola outbreak in West Africa. In this paper, we introduce vaccination of the 
susceptible population with the aim of controlling the spread of the disease 
and analyse two optimal control problems related with the transmission of Ebola 
disease with vaccination. Firstly, we consider the case where the total number 
of available vaccines in a fixed period of time is limited. Secondly, we analyse 
the situation where there is a limited supply of vaccines at each instant 
of time for a fixed interval of time. The optimal control problems have been 
solved analytically. Finally, we have performed a number of numerical simulations 
in order to compare the models with vaccination and the model without vaccination, 
which has recently been shown to fit the real data. Three 
vaccination scenarios have been considered for our numerical simulations, namely: 
unlimited supply of vaccines; limited total number of vaccines; and limited 
supply of vaccines at each instant of time.
\end{abstract}

\maketitle


\section{Introduction}
\label{sec:1}

Ebola is a lethal virus for humans that is currently under strong research 
due to the recent outbreak in West Africa and its socioeconomic impact 
(see, e.g., \cite{atangana2014,Gire20141369,Hayden2014,Kaushik2016254,%
HYG:996392,lekone,crimean,rivers2014,Walsh} and references therein). 
World Health Organization (WHO) has declared Ebola virus disease epidemic 
as a public health emergency of international concern with severe 
global economic burden. At fatal Ebola infection stage, patients usually die 
before the antibody response. Mainly after the 2014 Ebola outbreak in West Africa, 
some attempts to obtain a vaccine for Ebola disease have been realized. 
According to the WHO, results in July 2015 from an interim analysis of the Guinea 
Phase III efficacy vaccine trial show that VSV-EBOV (Merck, Sharp \& Dohme) 
is highly effective against Ebola \cite{Althaus}.

Since 2014, different mathematical models to analyze the spread of the 2014 
Ebola outbreak have been presented (see, e.g., 
\cite{MR3394468,Area:in:press,MyID:335,MR3349757,MyID:331} and references therein). 
In these models the populations under study are divided into compartments, 
and the rates of transfer between compartments are expressed mathematically 
as derivatives with respect to time of the size of the compartments. In a 
recent work \cite{Area:in:press}, a system of eight nonlinear (fractional) 
differential equations for a population divided into eight mutually exclusive 
groups was considered: susceptible, exposed, infected, hospitalized, asymptomatic 
but still infectious, dead but not buried, died, and completely recovered. 
By comparing the numerical results of this mathematical model and the real data 
provided by WHO, the difference in the period of 438 days analyzed is about 
7 cases per day. Note that in the day 438 after the beginning of the outbreak, 
the number of confirmed cases is 15018.

There exist different models for the spreading of Ebola, beginning with 
the simplest SIR and SEIR models \cite{Althaus,Chowell,ChowellN} and later 
more complex but also more realistic models have been considered 
\cite{Area:in:press,HYG:996392,Ngwa}. In \cite{lekone}, a stochastic discrete-time 
Susceptible-Exposed-Infectious-Recovered (SEIR) model for infectious diseases 
is developed with the aim of estimating parameters from daily incidence and mortality 
time series for an outbreak of Ebola in the Democratic Republic of Congo in 1995. 
In \cite{HYG:996392}, the authors use data from two epidemics (in Democratic Republic 
of Congo in 1995 and in Uganda in 2000) and built a SEIHFR 
(Susceptible-Exposed-Infectious-Hospitalized-F(dead but not yet buried)-Removed) 
mathematical model for the spread of Ebola haemorrhagic fever epidemics taking 
into account transmission in different epidemiological settings (in the community, 
in the hospital, during burial ceremonies). In \cite{atangana2014}, the authors 
propose a SIRD (Susceptible-Infectious-Recovered-Dead) mathematical model 
using classical and beta derivatives. In this model, the class of susceptible 
individuals does not consider new born or immigration. The study shows that, 
for small portion of infected individuals, the whole country could die out 
in a very short period of time in case there is no good prevention. 
In \cite{MR3394468}, a fractional order SEIR Ebola epidemic model 
is proposed and the authors show that the model gives 
a good approximation to real data published by WHO, 
starting from March 27th, 2014.    

Optimal control is a mathematical theory that emerged after 
the Second World War with the formulation of the celebrated 
Pontryagin maximum principle, responding to practical needs 
of engineering, particularly in the field of aeronautics 
and flight dynamics \cite{Pontryagin:1962}. In the last decade, 
optimal control has been largely applied to biomedicine, namely 
to models of cancer chemotherapy (see, e.g., 
\cite{Ledzewicz:cancer:SIAM:2007}), and recently to epidemiological models 
\cite{OC:HIV:PLoSCompBio:2015,OC:HepatiticC:PLoSOne:2011,SilvaTorres:TBHIV:2015}. 

In \cite{MR3349757}, the authors present a comparison between SIR and SEIR 
mathematical models used in the description of the Ebola virus propagation. 
They applied optimal control techniques in order to understand how the spread 
of the virus may be controlled, e.g., through education campaigns, 
immunization or isolation. In \cite{Dure}, the authors introduce 
a deterministic SEIR type model with additional hospitalization, 
quarantine and vaccination components in order to understand 
the disease dynamics. Optimal control strategies, both in the case 
of hospitalization (with and without quarantine) and vaccination, 
are used to predict the possible future outcome in terms 
of resource utilization for disease control and the effectiveness 
of vaccination on sick populations. Both in \cite{Dure} and \cite{MR3349757}, 
the authors study optimal control problems with $L^{2}$ cost 
functionals without any state or control constraints.
Here, we modify the model analyzed in \cite{Area:in:press} 
in order to consider optimal control problems with vaccination constraints. 
More precisely, we introduce an extra variable for the number of vaccines used, 
and we compare the hypothetical results if the vaccine were available 
at the beginning of the outbreak with the results of the model without vaccines. 
Firstly, we consider an optimal control problem with an end-point state constraint, 
that is, the total number of available vaccines, in a fixed period of time, 
is limited. Secondly, we analyze an optimal control problem with a mixed 
state constraint, in which there is a limited supply of vaccines at each instant 
of time for a fixed interval of time. Both optimal control problems have been 
analytically solved. Moreover, we have performed a number of numerical simulations 
in three different scenarios: unlimited supply of vaccines; limited total number 
of vaccines to be used; and limited supply of vaccines at each instant of time. 
From the results obtained in the first two cases, when there is no limit 
in the supply of vaccines or when the total number of vaccines used is limited, 
the optimal vaccination strategy implies a vaccination of 100\% of the susceptible 
population in a very short period of time (smaller than one day). In practice, 
this is a very difficult task because limitations in the number of vaccines 
and also in the number of humanitarian and medical teams in the affected regions
are common. In this direction, the third analyzed case is extremely important 
since we consider a limited supply of vaccines at each instant of time.

The paper is organized as follows. In Section~\ref{sec:2}, we recall 
a mathematical model for Ebola virus. In Section~\ref{sec:3}, the 
introduction of effective vaccination for Ebola virus is modeled. 
An optimal control problem with an end-point state constraint is formulated and solved analytically
in Section~\ref{sec:4}, which models the case where the total number of available vaccines 
in a fixed period of time is limited. In Section~\ref{sec:5}, the limited supply of vaccines 
at each instant of time for a fixed interval of time is mathematically translated into 
an optimal control problem with a mixed state control constraint. A closed form of 
the unique optimal control is given. In Section~\ref{sec:6}, we solve numerically 
the optimal control problems proposed in Sections~\ref{sec:4} and \ref{sec:5}. 
Finally, we end with Section~\ref{sec:7} of discussion of the results.


\section{Initial mathematical model for Ebola}
\label{sec:2}

\begin{figure}[ht!]
\centering
\includegraphics[width=0.75\textwidth]{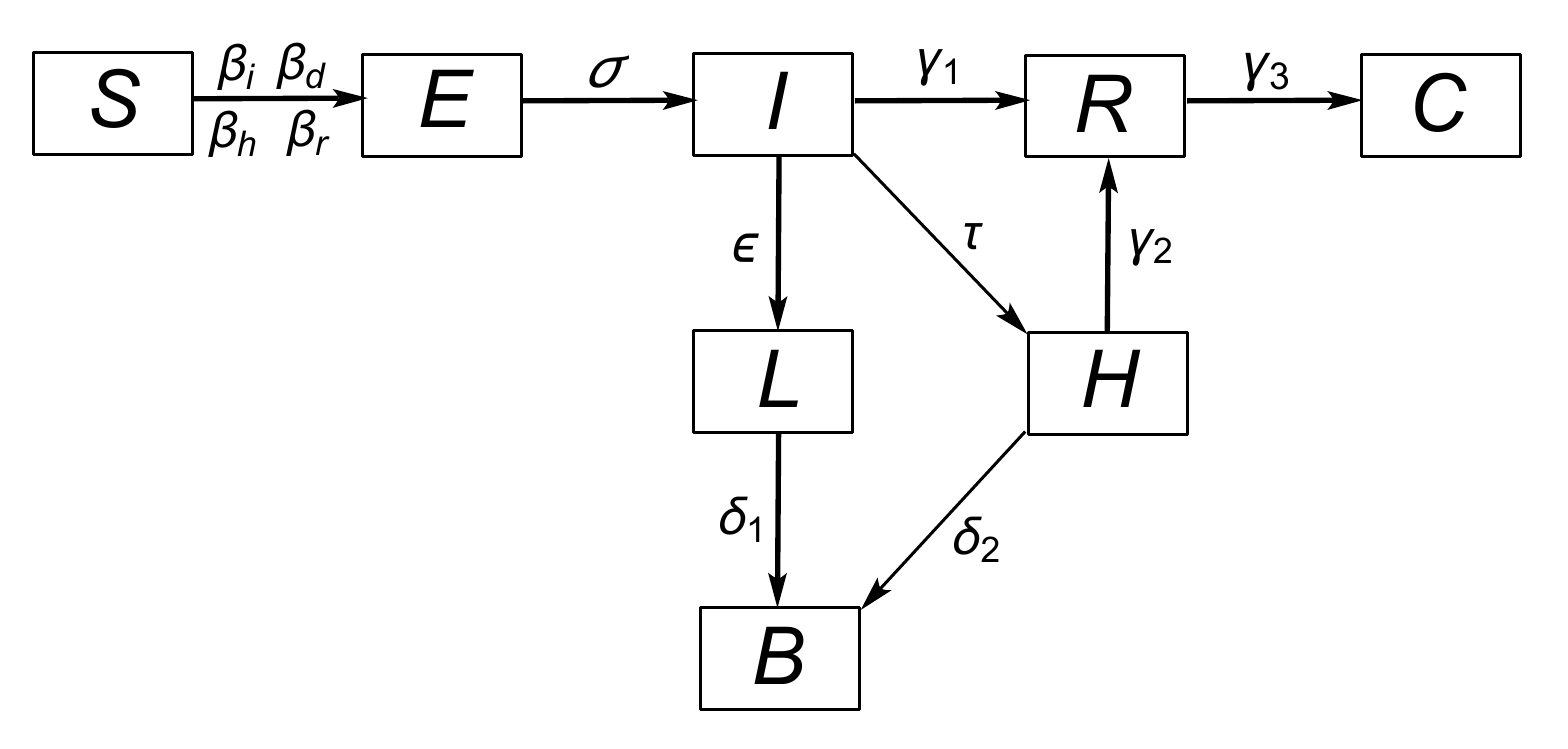}
\caption{Flowchart presentation of the compartmental model \eqref{eq:model} for Ebola.}
\label{figure:1}
\end{figure}

The total population $N$ under study is subdivided into eight mutually exclusive 
groups: susceptible ($S$), exposed ($E$), infected ($I$), hospitalized ($H$), 
asymptomatic but still infectious ($R$), dead but not buried ($D$), buried ($B$), 
and completely recovered ($C$). This model is adapted from \cite{ebola:models} 
and analyzed in \cite{Area:in:press}, where the birth and death rate are assumed
to be equal and are denoted by $\mu$, and the contact rate of susceptible 
individuals with infective, dead, hospitalized and asymptomatic individuals 
are denoted by $\beta_i$, $\beta_d$, $\beta_h$ and $\beta_r$, respectively. 
Exposed individuals become infectious at a rate $\sigma$. The per capita rate 
of progression of individuals from the infectious class to the asymptomatic
and hospitalized classes are denoted by $\gamma_1$ and $\tau$, respectively. 
Individuals in the dead class progress to the buried class at a rate $\delta_1$. 
Hospitalized individuals progress to the buried class and to the asymptomatic 
class at rates $\delta_2$ and $\gamma_2$, respectively. Asymptomatic individuals 
become completely recovered at a rate $\gamma_3$. Infectious individuals progress 
to the dead class at a fatality rate $\epsilon$. Dead and buried bodies are 
incinerated at a rate $\xi$. We assume that the total population, 
$N = S + E + I + R + H + D + B + C$, is constant, that is, the birth 
and death rates, both denoted by $\mu$, are equal to the incineration rate $\xi$. 
The model is mathematically described by the following system of 
eight nonlinear ordinary differential equations:
\begin{equation}
\label{eq:model}
\begin{cases}
\displaystyle \frac{dS}{dt} = \mu N - \frac{\beta_i}{N} S I - \frac{\beta_h}{N} S H
- \frac{\beta_d}{N} S D - \frac{\beta_r}{N} S R - \mu S,\\[0.2cm]
\displaystyle \frac{dE}{dt} = \frac{\beta_i}{N} S I + \frac{\beta_h}{N} S H
+ \frac{\beta_d}{N} S D + \frac{\beta_r}{N} S R - \sigma E - \mu E,\\[0.2cm]
\displaystyle \frac{dI}{dt} = \sigma E - (\gamma_1 + \epsilon + \tau + \mu)I,\\[0.2cm]
\displaystyle \frac{dR}{dt} = \gamma_1 I + \gamma_2 H - (\gamma_3 + \mu) R,\\[0.2cm]
\displaystyle \frac{dD}{dt} = \epsilon I - \delta_1 D - \xi D,\\[0.2cm]
\displaystyle \frac{dH}{dt} = \tau I - (\gamma_2 + \delta_2 + \mu) H,\\[0.2cm]
\displaystyle \frac{dB}{dt} = \delta_1 D + \delta_2 H - \xi B \\[0.2cm]
\displaystyle \frac{dC}{dt} = \gamma_3 R - \mu C.
\end{cases}
\end{equation}
In Fig.~\ref{figure:1}, we give a flowchart presentation of model 
\eqref{eq:model}. In this flowchart, we identify the compartmental classes 
as well as the parameters appearing in the model. Moreover, the values 
of the parameters are given in Table~\ref{table:parameters}.

The basic reproduction number (that is, the number of cases one case generates 
on average over the course of its infectious period, in an otherwise 
uninfected population) of model \eqref{eq:model} can be computed using the 
associated next-generation matrix method \cite{diekmann}. It is obtained 
as the spectral radius of the following matrix, known as the next-generation-matrix:
\begin{equation*}
FV^{-1} =  \begin{pmatrix}
A_{11} & A_{12} & A_{13} & A_{14} & A_{15}   \\
0 &0 & 0 & 0 &0   \\
0 &0 & 0 & 0 &0   \\
0 &0 & 0 & 0 &0   \\
0 &0 & 0 & 0 &0   \\
0 &0 & 0 & 0 &0   \\
\end{pmatrix},
\end{equation*}
where 
\[
A_{11} = \frac{a_3 \beta_i \sigma a_1 a_4 + a_3 \beta_r  \sigma 
\left(a_4\gamma_1 +\tau \gamma_2\right)+\beta_d  \epsilon \sigma 
a_1 a_4 +a_3 \beta_h  \tau \sigma a_1}{ a_1 a_2 a_3 a_4 a_5} ,
\]
\[ 
A_{12}= \frac{\beta_r }{a_1}+\frac{\beta_r (a_4\gamma_1 + \tau \gamma_2)}{a_1 a_2 a_3} 
+ \frac{\beta_d  \epsilon}{a_2 a_3}+\frac{\beta_h  \tau}{a_2 a_4}, 
\]
\[ 
A_{13}= \frac{\beta_r }{a_1}, \quad
A_{14}=\frac{\beta_d  }{a_3}, \quad 
A_{15}=\frac{\beta_r \gamma_2}{a_1 a_4}+ \frac{\beta_h }{a_4}
\]
with 
\[
a_1 = \gamma_2 + \delta_2 + \mu,
\quad a_2=  \gamma_3 + \mu, 
\quad a_3 = \delta_1 + \xi, 
\quad a_4= \sigma + \mu, 
\quad \text{and} \quad  a_5 =\gamma_1 + \epsilon + \tau + \mu. 
\]
Therefore, the basic reproduction number $R_0$ is given by
\begin{multline*}
R_0 = \frac{\sigma}{a_1 a_2 a_3 a_4 a_5} \left[
a_3 \left(\beta_{i} (\gamma_3 a_1
+ \mu (\gamma_2 + \delta_2) +{\mu}^{2})
+\beta_{r} (\gamma_1 a_1
+ \gamma_2 \tau)
+\beta_{h} \tau a_2 \right)
\right.  \\
\left. +\beta_{d} \epsilon \left((\gamma_3 a_1
+ \mu (\gamma_2 + \delta_2))
+{\mu}^{2} \right)
\right].
\end{multline*}
As it is well-known, if the basic reproduction number $R_{0}<1$, 
then the infection will stop in the long run; but if $R_{0}>1$, 
then the infection will spread in population.

\begin{table}
\begin{adjustwidth}{-0.5in}{0in}
\begin{tabular}{|l|l|l|}
\hline
{\small{Symbol}} & {\small{Description}}  & {\small{Value}} \\
\hline
{\small{$\sigma$}} & {\small{per capita rate at which exposed
individuals become infectious}}  & {\small{$1/11.4$}}\\
{\small{$\mu$}} & {\small{death rate}}  & {\small{$14/1000$}}\\
{\small{$\beta_i$}}
& {\small{contact rate of infective and susceptible individuals}}
& {\small{$0.14$}}\\
{\small{$\beta_d $}}
& {\small{contact rate of infective and dead individuals}}
& {\small{$0.40 $}}\\
{\small{$\beta_h$}}
& {\small{contact rate of infective and hospitalized individuals}}
& {\small{$0.29$}}\\
{\small{$\beta_r$}}
& {\small{contact rate of infective and asymptomatic individuals}}
& {\small{$0.185$}}\\
{\small{$\gamma_1$}}
& {\small{per capita rate of progression of individuals from the infectious class}}
&  \\
& {\small{to the asymptomatic class}} & {\small{$1/10$}}\\
{\small{$\epsilon$}} & {\small{fatality rate }} & {\small{$1/9.6$}}\\
{\small{$\delta_1$}}
& {\small{per capita rate of progression of individuals from the dead class}} & \\
& {\small{ to the buried class}} & {\small{$1/2$}}\\
{\small{$\delta_2$}}
& {\small{per capita rate of progression of individuals from the hospitalized class}}
& \\
& {\small{to the buried class}} & {\small{$1/4.6$}}\\
{\small{$\gamma_2$}}
& {\small{per capita rate of progression of individuals from the hospitalized class }}
& \\
& {\small{to the asymptomatic class}} & {\small{$1/5$}}\\
{\small{$\tau$}}
& {\small{per capita rate of progression of individuals from the infectious class}}
& \\
& {\small{to the hospitalized class}} & {\small{$1/5$}}\\
{\small{$\gamma_3$}}
& {\small{per capita rate of progression of individuals from the asymptomatic class}} & \\
& {\small{to the completely recovered class}} & {\small{$1/30$}}\\
{\small{$\xi$}} & {\small{incineration rate}} & {\small{$14/1000$}}\\
\hline
\end{tabular}
\caption{Parameter values for model \eqref{eq:model}, corresponding
to a basic reproduction number $R_0 = 2.287$. The values of the parameters 
come from \cite{bwaka,dowell,khan,HYG:996392,ndambi,rivers2014,rowe}.}
\label{table:parameters}
\end{adjustwidth}
\end{table}

In this section, we have recalled a model for describing 
the Ebola virus transmission. Now we want to address the question 
about how to introduce vaccination as a prevention measure. 
This is analyzed in the next section.


\section{Mathematical model for Ebola with vaccination}
\label{sec:3}

We now introduce vaccination of the susceptible population with the aim 
of controlling the spread of the disease. We assume that the vaccine 
is effective so that all vaccinated susceptible individuals become completely 
recovered (see, e.g., \cite{SEIR:Rosario:2014,NeilanLenhart2010} 
for vaccination in a SEIR model that corresponds to a system of four 
nonlinear ordinary differential equations). Let us introduce in model 
\eqref{eq:model} a control function $u(t)$, which represents the percentage 
of susceptible individuals being vaccinated at each instant 
of time $t$ with $t \in [0, t_f]$. Unfortunately, in many situations, the number 
of available vaccines does not fulfill the necessities in order to eradicate 
the disease. In this paper, we consider limitations on the total number 
of vaccines during a fixed interval of time $[0, t_f]$ and on the number 
of available vaccines at each instant of time $t$ with $t \in [0, t_f]$. 
In order to translate this real situation mathematically, we introduce an extra
variable $W$ that denotes the number of vaccines used:
\begin{equation*}
\frac{d W}{dt}(t) = u(t) S(t),
\quad \text{subject to the initial condition } W(0) = 0.
\end{equation*}
Hence, the model with vaccination is given by the following system of 
nine nonlinear ordinary differential equations:
\begin{equation}
\label{eq:model:vaccine}
\begin{cases}
\displaystyle \frac{dS}{dt}(t) = \mu N - \frac{\beta_i}{N} S(t) I(t) 
- \frac{\beta_h}{N} S(t) H(t) - \frac{\beta_d}{N} S(t) D(t)\\[0.2cm] 
\displaystyle \hspace*{1cm}- \frac{\beta_r}{N} S(t) R(t) - \mu S(t) - S(t) u(t),\\[3mm]
\displaystyle \frac{dE}{dt}(t) 
= \frac{\beta_i}{N} S(t) I(t) + \frac{\beta_h}{N} S(t) H(t)
+ \frac{\beta_d}{N} S(t) D(t)\\[3mm]
\displaystyle  \hspace*{1cm} + \frac{\beta_r}{N} S(t) R(t)  - \sigma E(t) - \mu E(t),\\[3mm]
\displaystyle \frac{dI}{dt}(t) = \sigma E(t) - (\gamma_1 + \epsilon + \tau + \mu)I(t),\\[3mm]
\displaystyle \frac{dR}{dt}(t) = \gamma_1 I(t) + \gamma_2 H(t) - (\gamma_3 + \mu) R(t),\\[3mm]
\displaystyle \frac{dD}{dt}(t) = \epsilon I(t) - (\delta_1 + \xi) D(t),\\[3mm]
\displaystyle \frac{dH}{dt}(t) = \tau I(t) - (\gamma_2 + \delta_2 + \mu) H(t),\\[3mm]
\displaystyle \frac{dB}{dt}(t) = \delta_1 D(t) + \delta_2 H(t) - \xi B(t),\\[3mm]
\displaystyle \frac{dC}{dt}(t) = \gamma_3 R(t) - \mu C(t) + S(t) u(t),\\[3mm]
\displaystyle \frac{dW}{dt}(t) = S(t) u(t).
\end{cases}
\end{equation}

In model \eqref{eq:model:vaccine}, the vaccination parameter is fixed. 
In the next section, we address the question of how to choose 
this parameter in an optimal way along time.


\section{Optimal control with an end-point state constraint}
\label{sec:4}

We start by considering the case where the total number of available vaccines, 
in a fixed period of time, is limited. We formulate and solve analytically such 
optimal control problem with end-point state constraint, 
which will be then solved numerically in Section~\ref{sec:6}.

We consider the model with vaccination \eqref{eq:model:vaccine} and formulate 
the optimal control problem with the aim to determine the vaccination strategy 
$u$ over a fixed interval of time $[0, t_f]$ that minimizes the cost functional
\begin{equation}
\label{eq:cost:vaccine:L2}
J(u) = \int_0^{t_f} \left[ w_1 I(t) + w_2 u^2(t) \right]  \, dt,
\end{equation}
where the constants $w_1$ and $w_2$ represent the weights associated with the 
number of infected individuals and on the cost  associated with the vaccination 
program, respectively. We assume that the control function $u$ takes values 
between 0 and 1. When $u(t)=0$, no susceptible individual is vaccinated at time 
$t$; if $u(t) =1$, then all susceptible individuals are vaccinated at $t$. 
Let $\vartheta$ denote the total amount of available vaccines in a fixed period 
of time $[0, t_f]$. This constraint is represented by
\begin{equation}
\label{totalconstraint}
W(t_f) \leq \vartheta.
\end{equation}
Let
\begin{equation*}
x(t) =(x_1(t), \ldots,  x_9(t))
=\left( S(t), E(t), I(t), R(t), D(t), H(t), B(t), C(t), W(t) \right) 
\in {\mathbb{R}}^9.
\end{equation*}
The optimal control problem consists to find the optimal trajectory $\tilde{x}$ 
associated with the optimal control $\tilde{u}$, satisfying the control system 
\eqref{eq:model:vaccine}, the initial conditions 
\begin{equation}
\label{initcond9}
x(0) = (18000, 0, 15, 0, 0, 0, 0, 0, 0)
\end{equation}
(see \cite{Area:in:press}), the constraint \eqref{totalconstraint}, 
and where the control $\tilde{u} \in \Omega$ minimizes the objective 
functional \eqref{eq:cost:vaccine:L2} with
\begin{equation}
\label{eq:admiss:control}
\Omega = \biggl\{ u(\cdot) \in L^{\infty}(0, t_f) \,
| \,  0 \leq u (t) \leq 1  \biggr\}.
\end{equation}
The existence of an optimal control $\tilde{u}$ and associated optimal trajectory 
$\tilde{x}$ comes from the convexity of the integrand of the cost function 
\eqref{eq:cost:vaccine:L2} with respect to the control $u$ and the Lipschitz 
property of the state system with respect to state variables $x$ (see, \textrm{e.g.}, 
\cite{Cesari_1983,Fleming_Rishel_1975} for existence results of optimal solutions). 
According to the Pontryagin Maximum Principle \cite{Pontryagin:1962}, 
if $\tilde{u} \in \Omega$ is optimal for the problem \eqref{eq:model:vaccine}, 
\eqref{eq:cost:vaccine:L2} with initial conditions \eqref{initcond9}
and fixed final time $t_f$, then there exists
$\lambda \in AC ([0, t_f]; \mathbb{R}^9)$,
$\lambda(t) = \left(\lambda_1(t), \ldots , \lambda_9(t)\right)$,
called the \emph{adjoint vector}, such that
\begin{equation*}
\dot{x} = \frac{\partial {\mathcal{H}}_{1}}{\partial \lambda}
\quad \quad \text{and} \quad \quad
\dot{\lambda} = -\frac{\partial {\mathcal{H}}_{1}}{\partial x},
\end{equation*}
where the Hamiltonian $\mathcal{H}_1$ is defined by
\begin{equation*}
\mathcal{H}_1(x, u, \lambda) =  w_1 x_3 + w_2 u^2(t)
+ \lambda \left(f(x) + A x + B xu  \right)
\end{equation*}
with
\begin{equation*}
f = \left(f_1 \, f_2 \, 0\, 0\, 0\, 0\, 0\, 0\, 0 \right),
\end{equation*}
\begin{equation*}
f_1 = \mu N - \frac{\beta_i}{N} S(t) I(t) - \frac{\beta_h}{N} S(t) H(t)
- \frac{\beta_d}{N} S(t) D(t) - \frac{\beta_r}{N} S(t) R(t),
\end{equation*}
\begin{equation*}
f_2 = \frac{\beta_i}{N} S(t) I(t) + \frac{\beta_h}{N} S(t) H(t)
+ \frac{\beta_d}{N} S(t) D(t) + \frac{\beta_r}{N} S(t) R(t),
\end{equation*}
\begin{equation*}
A =  \begin{pmatrix}
- \mu & 0 & 0 & 0 &0  & 0 & 0 & 0 \\
0 & - \sigma - \mu  & 0 & 0 &0  & 0 & 0 & 0 \\
0 & \sigma & \varLambda & 0 &0  & 0 & 0 & 0 \\
0 & 0 & \gamma_1 & - (\gamma_3 + \mu) &0  & \gamma_2 & 0 & 0 \\
0 & 0 & \epsilon & 0 & - (\delta_1 + \xi)  & 0 & 0 & 0 \\
0 & 0 & \tau & 0 &0  & - (\gamma_2 + \delta_2 + \mu) & 0 & 0 \\
0 & 0 & 0 & 0 & \delta_1  &  \delta_2 & - \xi & 0 \\
0 & 0 & 0 & \gamma_3 &0  & 0 & 0 & - \mu \\
0 & 0 & 0 & 0 &0  & 0 & 0 & 0 \\
\end{pmatrix}
\end{equation*}
with $\varLambda = - (\gamma_1 + \epsilon + \tau + \mu)$,
and
\begin{equation*}
B = \left(b \, Z\right),
\end{equation*}
where $b = \left(-1 \, 0 \, 0\, 0\, 0\, 0\,  0\, 1\, 1\right)^T$
and $Z = 0$ with $0$ the $8 \times 9$ null matrix.
The minimization condition
\begin{equation}
\label{maxcondPMP}
\mathcal{H}_1\left(\tilde{x}(t), \tilde{u}(t), \tilde{\lambda}(t)\right)
= \min_{u \in \Omega}
\mathcal{H}_1\left(\tilde{x}(t), u, \tilde{\lambda}(t)\right)
\end{equation}
holds almost everywhere on $[0, t_f]$. Moreover, transversality conditions
$\lambda_i(t_f) = 0$, $i =1,\ldots, 8$, hold. Solving the minimality condition
\eqref{maxcondPMP} on the interior of the set of admissible controls 
$\Omega$ gives 
\begin{equation*}
\tilde{u}(t) = -\frac{1}{2}\frac{\left(-\tilde{\lambda}_1(t) + \tilde{\lambda}_8(t)
+ \tilde{\lambda}_9(t) \right)\tilde{x}_1(t)}{\omega_2},
\end{equation*}
where the adjoint functions satisfy 
\begin{equation}
\label{adjointsystem}
\begin{cases}
\displaystyle \dot{\tilde{\lambda}}_1 
= -\tilde{\lambda}_1 \, \left( -\frac{\beta_i}{N}\tilde{x}_3
- \frac{\beta_h}{N}\tilde{x}_6-\frac{ \beta_d}{N}\tilde{x}_5
- \frac{ \beta_r}{N}\tilde{x}_4 -\mu - \tilde{u} \right)\\ 
\displaystyle \hspace*{0.9cm}
- \tilde{\lambda}_2 \left(\frac{\beta_{i}}{N}\tilde{x}_3
+ \frac{\beta_{h}}{N}\tilde{x}_6 +\frac{\beta_{d}}{N}\tilde{x}_5
+\frac{\beta_{r}}{N}\tilde{x}_4 \right) - \tilde{\lambda}_8 \tilde{u} 
- \tilde{\lambda}_9 \tilde{u}, \\[0.2 cm]
\dot{\tilde{\lambda}}_2 = -\tilde{\lambda}_2 \, \left( -\sigma-\mu \right)
- \tilde{\lambda}_3 \,\sigma, \\[0.2cm]
\displaystyle \dot{\tilde{\lambda}}_3 
= -\omega_1 + \tilde{\lambda}_1\frac{\beta_{i}}{N}\tilde{x}_1
-\tilde{\lambda}_2\frac{\beta_{i}}{N}\tilde{x}_1 - \tilde{\lambda}_3 \,
\left( - \gamma_1 - \epsilon-\tau-\mu \right) - \tilde{\lambda}_4 \gamma_1
- \tilde{\lambda}_5 \epsilon- \tilde{\lambda}_6 \tau, \\[0.2cm]
\displaystyle \dot{\tilde{\lambda}}_4 
=  \tilde{\lambda}_1\frac{\beta_{r}}{N}\tilde{x}_1
- \tilde{\lambda}_2 \frac{\beta_{r}}{N}\tilde{x}_1 -\tilde{\lambda}_4
\left( -\gamma_3 -\mu \right) - \tilde{\lambda}_8 \gamma_3, \\[0.2cm]
\displaystyle \dot{\tilde{\lambda}}_5 
=  \tilde{\lambda}_1 \frac{\beta_{d}}{N}\tilde{x}_1
- \tilde{\lambda}_2 \frac{\beta_{d}}{N}\tilde{x}_1- \tilde{\lambda}_5
\left( -\delta_1 -\xi \right) - \tilde{\lambda}_7 \delta_1,\\[0.2cm]
\displaystyle \dot{\tilde{\lambda}}_6 
= \tilde{\lambda}_1\frac{\beta_{h}}{N}\tilde{x}_1
- \tilde{\lambda}_2\frac{\beta_{h}}{N}\tilde{x}_1 - \tilde{\lambda}_4 \gamma_2
- \tilde{\lambda}_6 \left(-\gamma_2 - \delta_2 -\mu\right)
- \tilde{\lambda}_7 \delta_2,\\[0.2cm]
\displaystyle \dot{\tilde{\lambda}}_7 = \tilde{\lambda}_7 \xi, \\[0.2cm]
\displaystyle \dot{\tilde{\lambda}}_8 = \tilde{\lambda}_8 \mu, \\[0.2cm]
\displaystyle \dot{\tilde{\lambda}}_9 =  0.
\end{cases}
\end{equation}
Since $W$ has initial and terminal conditions, the adjoint function $\lambda_9$,
 associated with the state variable $W$, has no transversality condition.
From \eqref{adjointsystem}, $\tilde{\lambda}_9 \equiv k$, where the constant $k$
must be such that the end point conditions
$W(0) = 0$ and $W(t_f) = \vartheta$ are satisfied. As the optimal control
$\tilde{u}$ can take values on the boundary of the control set $[0, 1]$,
the optimal control $\tilde{u}$ must satisfy
\begin{equation}
\label{optcontrol}
\tilde{u}(t) = \min \left\{1,
\max \left\{0,  \frac{1}{2}\frac{\left(\tilde{\lambda}_1(t)
- \tilde{\lambda}_8(t) - \tilde{\lambda}_9(t) \right)
\tilde{x}_1(t)}{\omega_2}  \right\} \right\}.
\end{equation}
The optimal control $\tilde{u}$ given by \eqref{optcontrol}
is unique due to the boundedness of the state and adjoint functions
and the Lipschitz property of systems \eqref{eq:model:vaccine}
and \eqref{adjointsystem}.

We would like to note that if we consider the optimal control 
problem without any restriction on the number of available vaccines, 
that is, to find the optimal solution
$(\tilde{x}, \tilde{u})$, with $\tilde{u} \in \Omega$, which minimizes
the cost functional \eqref{eq:cost:vaccine:L2} subject to the control system
\eqref{eq:model:vaccine}, initial conditions \eqref{initcond9}, and free final
conditions $(x_1(t_f), \ldots, x_9(t_f))$, then the adjoint functions
$(\lambda_1, \ldots, \lambda_9)$ must satisfy transversality conditions
$\lambda_i(t_f) = 0$, $i =1,\ldots, 9$, and, since $\tilde{\lambda}_9 = 0$,
the optimal control is given by
\begin{equation*}
\tilde{u}(t) = \min \left\{1, \max \left\{0, \frac{1}{2}\frac{\left(\tilde{\lambda}_1(t)
- \tilde{\lambda}_8(t) \right)\tilde{x}_1(t)}{\omega_2} \right\} \right\}.
\end{equation*}

In a concrete situation, the number of available vaccines is always limited. 
Therefore, it is also important to study the optimal control problem with 
such kind of constraints. This is done in Section~\ref{sec:5}. Both problems, 
with and without constraints, are  numerically solved in Section~\ref{sec:6}.


\section{Optimal control with a mixed state control constraint}
\label{sec:5}

A particularly challenging situation in vaccination programs happens when 
there is a limited supply of vaccines at each instant of time for a fixed 
interval of time $[0, t_f]$. In order to study this health public problem, 
from the optimal point of view, we formulate an optimal control problem with 
a mixed state control constraint (see, e.g., \cite{SEIR:Rosario:2014}).
The cost functional \eqref{eq:cost:vaccine:L2} remains \eqref{eq:cost:vaccine:L2}, 
the one considered in previous section, as well as the set of admissible controls 
$\Omega$ \eqref{eq:admiss:control}. The end point state constraint 
\eqref{totalconstraint} is replaced by the following mixed state control constraint:
\begin{equation*}
\label{eq:constraint}
S(t) u(t) \leq \vartheta \, , \quad \vartheta \geq 0 \, , \, \,
\text{for almost all} \, \, t \in [0, t_f],
\end{equation*}
which should be satisfied at almost every instant of time during the whole 
vaccination program. Analogously to \cite{SEIR:Rosario:2014}, we observe 
that in our optimal control problem the differential equation
\begin{equation*}
\frac{dW}{dt}(t) = S(t) u(t)
\end{equation*}
does not appear neither in the cost and in any other differential equation,  
nor in the mixed state control constraint. Thus, in this section, the control 
system does not include the last equation and $x$ is used to denote
\begin{equation*}
x(t) =(x_1(t), \ldots,  x_8(t))
=\left( S(t), E(t), I(t), R(t), D(t), H(t), B(t), C(t) \right) 
\in {\mathbb{R}}^8.
\end{equation*}
Let us consider the initial conditions \eqref{initcond9}. 
The control system can be rewritten in the following way:
\begin{equation*}
\frac{dx(t)}{dt} = f(x(t)) + A x(t)+ B x(t) u(t),
\end{equation*}
with
\begin{equation*}
A = \begin{pmatrix}
- \mu & 0 & 0 & 0 &0  & 0 & 0 & 0 \\
0 & - \sigma - \mu  & 0 & 0 &0  & 0 & 0 & 0 \\
0 & \sigma & \varLambda & 0 &0  & 0 & 0 & 0 \\
0 & 0 & \gamma_1 & - (\gamma_3 + \mu) &0  & \gamma_2 & 0 & 0 \\
0 & 0 & \epsilon & 0 & - (\delta_1 + \xi)  & 0 & 0 & 0 \\
0 & 0 & \tau & 0 &0  & - (\gamma_2 + \delta_2 + \mu) & 0 & 0 \\
0 & 0 & 0 & 0 & \delta_1  &  \delta_2 & - \xi & 0 \\
0 & 0 & 0 & \gamma_3 &0  & 0 & 0 & - \mu \\
\end{pmatrix},
\end{equation*}
$\varLambda = - (\gamma_1 + \epsilon + \tau + \mu)$, 
\begin{equation*}
B = \left(b \, Z  \right),
\end{equation*}
where $b = \left(-1 \, 0 \, 0\, 0\, 0\, 0\,  0\, 1 \right)^T$ 
and $Z = 0$ with $0$ the $7 \times 8$ null matrix, and
\begin{equation*}
f =  \left(f_1 \, f_2 \, 0\, 0\, 0\, 0\, 0\, 0 \right),
\end{equation*}
with
\begin{equation*}
f_1 = \mu N - \frac{\beta_i}{N} S(t) I(t) - \frac{\beta_h}{N} S(t) H(t)
- \frac{\beta_d}{N} S(t) D(t) - \frac{\beta_r}{N} S(t) R(t)
\end{equation*}
and
\begin{equation*}
f_2 = \frac{\beta_i}{N} S(t) I(t) + \frac{\beta_h}{N} S(t) H(t)
+ \frac{\beta_d}{N} S(t) D(t) + \frac{\beta_r}{N} S(t) R(t).
\end{equation*}
It follows from Theorem~23.11 in \cite{livro:Clarke:2013} that
problem \eqref{eq:model:vaccine}--\eqref{eq:admiss:control}
has a solution (see also \cite{SEIR:Rosario:2014}).
Let $(\tilde{x}, \tilde{u})$ denote such solution.
To determine it, we apply the Pontryagin Maximum Principle
(see, e.g., Theorem~7.1 in \cite{SIAM:Clarke:Rosario}): there exists
multipliers $\lambda_0 \leq 0$, $\lambda \in AC([0, t_f]; {\mathbb{R}}^8)$,
and $\psi \in L^1([0, t_f]; {\mathbb{R}})$, such that
\begin{itemize}
\item $\min \{ | \lambda(t) | \, : \, t \in [0, t_f] \}
> \lambda_0$ (nontriviality condition);

\item $\frac{d \lambda (t)}{dt}
= - \frac{\partial \mathcal{H}_2}{\partial x}(\tilde{x}(t), 
\tilde{u}(t), \lambda_0, \lambda(t), \psi(t))$ (adjoint system);

\item $\lambda(t) B \tilde{x}(t) + \psi(t) \tilde{x}_1(t)
+ \lambda_0 w_2 \tilde{u}^2(t) \in \mathcal{N}_{[0,1]}(\tilde{u}(t))$ a.e. and
$$
\mathcal{H}_2(\tilde{x}(t), \tilde{u}(t), \lambda_0, \lambda(t), \psi(t))
\leq \mathcal{H}_2(\tilde{x}(t), v, \lambda_0, \lambda(t), \psi(t)), \,
\forall v \in [0, 1] \, :\,  \tilde{x}_1(t)v \leq \vartheta
$$
(minimality condition);

\item $\psi(t) (\tilde{x}_1(t) \tilde{u}(t) - \vartheta) = 0$
and $\psi(t) \leq 0$ a.e.;

\item $\lambda(t_f) = (0, \ldots, 0)$ (transversality conditions);
\end{itemize}
where the Hamiltonian $\mathcal{H}_2$ for problem
\eqref{eq:model:vaccine}--\eqref{eq:admiss:control} is defined by
\begin{equation*}
\mathcal{H}_2(x, u, \lambda_0, \lambda, \psi)
= \lambda_0 \left( w_1 x_3 + w_2 u^2 \right)
+ \lambda \left( f(x) + A x + B x u \right) + \psi (S u - \vartheta),
\end{equation*}
and
$\mathcal{N}_{[0,1]}(\tilde{u}(t))$ stands for the normal cone from convex
analysis to $[0, 1]$ at the optimal control $\tilde{u}(t)$
(see, e.g., \cite{livro:Clarke:2013}). The optimal solution
$(\tilde{x}, \tilde{u})$ is normal (see \cite{SEIR:Rosario:2014} for details),
so we can choose $\lambda_0 = 1$. Analogously to previous section, we obtain
a closed form of the unique optimal control $\tilde{u}$:
\begin{equation*}
\tilde{u}(t) = \min \left\{1, \max \left\{0,  \frac{1}{2}\frac{\left(\tilde{\lambda}_1(t)
- \tilde{\lambda}_8(t) - \psi(t) \right)\tilde{x}_1(t)}{\omega_2} \right\} \right\}.
\end{equation*}

The theoretical results obtained in Sections~\ref{sec:4} and \ref{sec:5} 
are illustrated numerically in the next section.


\section{Numerical simulations}
\label{sec:6}

\begin{figure}[!htb]
\centering
\subfloat[\footnotesize{Real data and solution of \eqref{eq:model}}]{
\label{Cum:Inf:data}\includegraphics[width=0.45\textwidth]{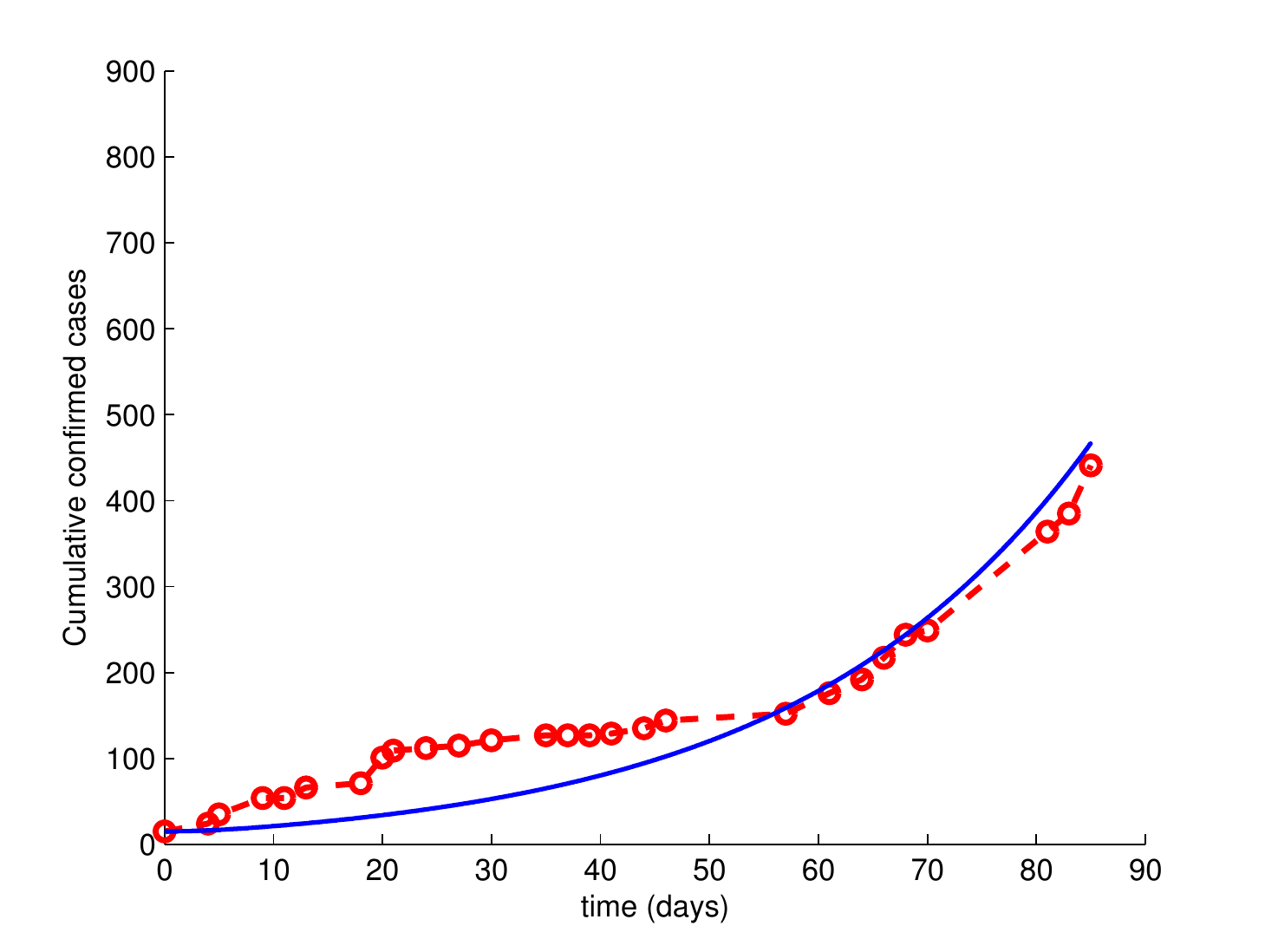}}
\subfloat[\footnotesize{Unlimited supply of vaccines}]{
\label{Cum:Inf:VacNoLim}\includegraphics[width=0.45\textwidth]{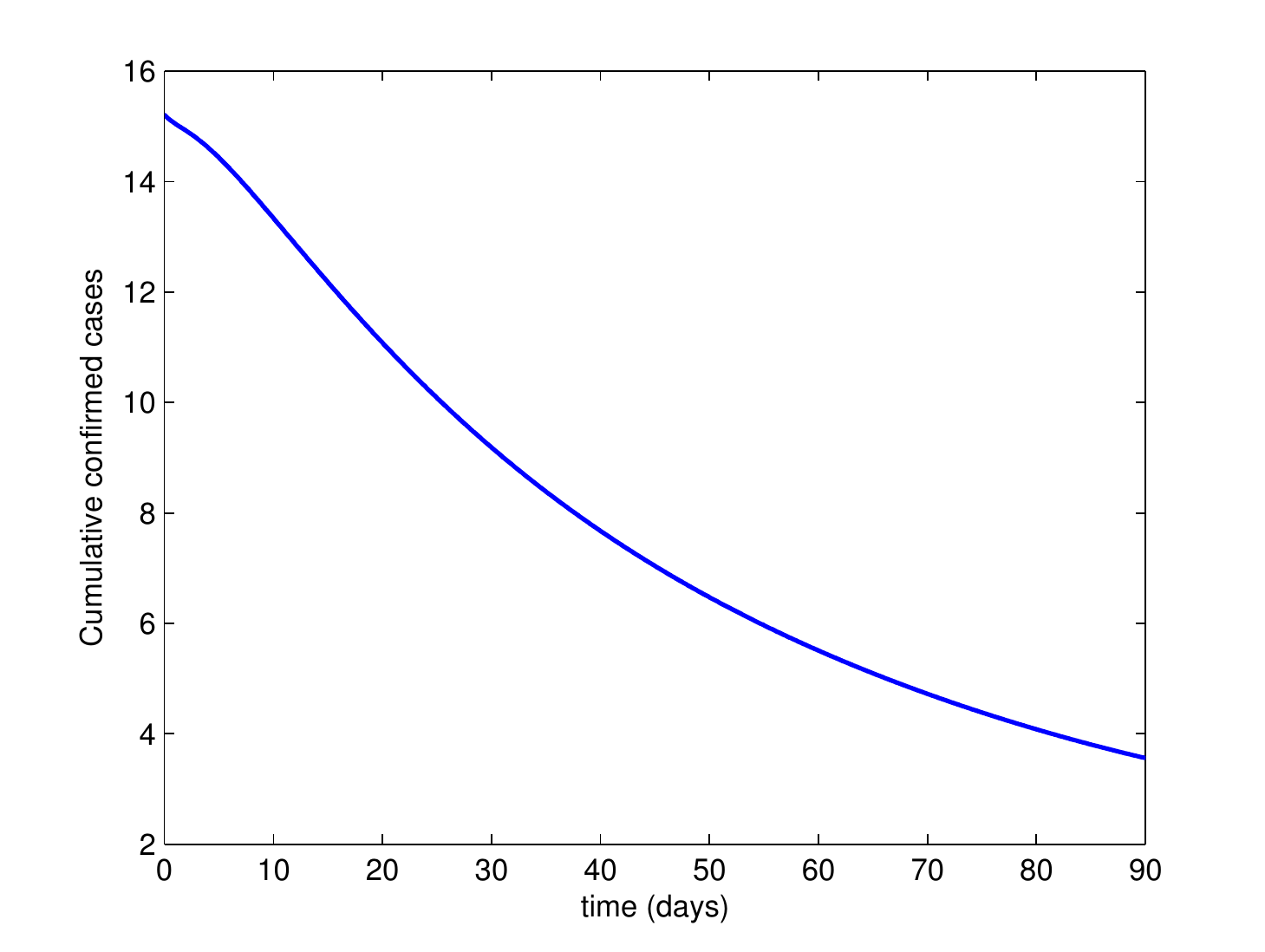}}
\caption{(a) Cumulative confirmed cases: in dashed circle line the real data
from WHO and in continuous line the values of
$I(t) + R(t) + D(t) + H(t) + B(t) + C(t) - \mu(N - S(t) - E(t))$
from \eqref{eq:model} with the parameter values from 
Table~\ref{table:parameters}. (b) Cumulative confirmed cases given in
\eqref{eq:model:vaccine}, when available an unlimited supply 
of vaccines, also with the parameter values 
from Table~\ref{table:parameters}.}
\label{fig:cum:Inf}
\end{figure}

We start the numerical simulations by considering an intervention of 90 days, 
initial conditions given in \eqref{initcond9}, and the evolution of cumulative 
confirmed cases based on the data from the World Health Organization (WHO), 
following all the reports of the disease in the three main affected countries 
of Western Africa of the 2014 Ebola outbreak, namely, Liberia, Guinea and Sierra 
Leone. The model \eqref{eq:model} with the parameter values from 
Table~\ref{table:parameters} fits the real data from WHO, see Fig.~\ref{Cum:Inf:data}. 
We would like to emphasize that we are just considering the initial period 
of spreading of the disease, in which the vaccination should be introduced. 
Looking to a longer period of time (as considered in \cite{Area:in:press}), 
then the model fits quite well the real data: in \cite[Figure 2]{Area:in:press} 
the $\ell_{2}$ norm of the difference between the real data and the prediction 
is 3181, which gives an error of less than 7.3 cases per day, 
as compared with about 15,000 cases at the end of the outbreak.

Assuming that, in the near future, an effective vaccine against 
the Ebola virus will be available, as expected by WHO 
by the end of 2015, we study three different scenarios, 
which illustrate limitations on the number of vaccines available 
and on the capacity of administration of the vaccines by the health 
care services and humanitarian teams working in the affected countries. 
The three vaccination scenarios are the following: unlimited supply of vaccines; 
limited total number of vaccines to be used; and limited supply 
of vaccines at each instant of time.


\subsection{Unlimited supply of vaccines}
\label{sec:6.1}

Assume $w_1 = w_2 = 1$.  If we administer from an unlimited supply 
of vaccines, then the number of total individuals who have 
an active infection $I(t) + R(t) + D(t) + H(t) + B(t) 
- \mu(N - S(t) - E(t) - C(t))$ during the 90 days 
is a decreasing function in time, and is equal to $3.56$
individuals at the final time (see Fig.~\ref{Cum:Inf:VacNoLim}).

\begin{figure}[!htb]
\centering
\includegraphics[width=0.75\textwidth]{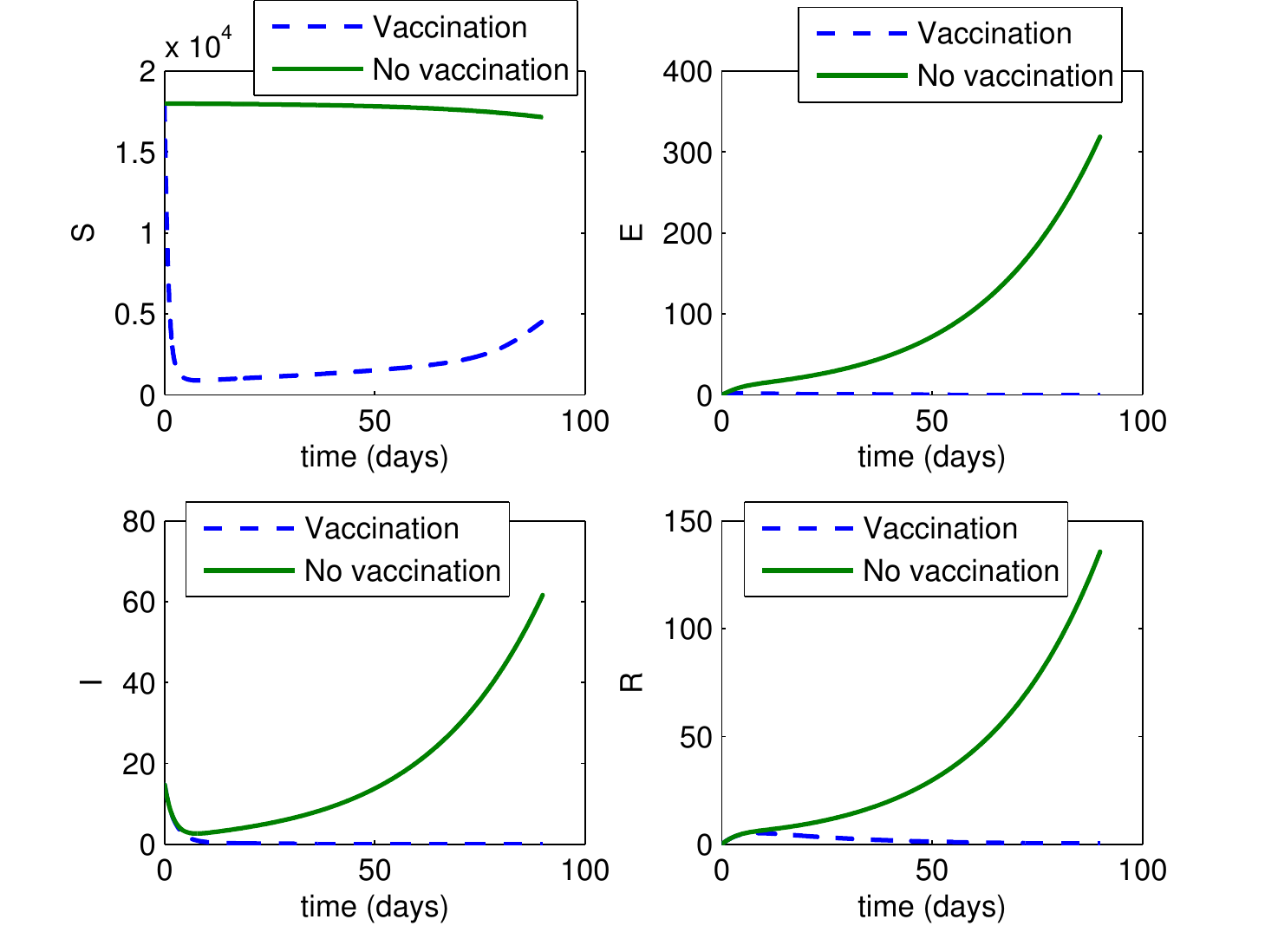}
\caption{Individuals $S(t)$, $E(t)$, $I(t)$ and $R(t)$.
In dashed line, the case of vaccination without limit
on the supply of vaccines; in continuous line,
the case with no vaccination with the parameter 
values from Table~\ref{table:parameters}.}
\label{fig:SEIR:vac:novac}
\end{figure}

\begin{figure}[!htb]
\centering
\includegraphics[width=0.75\textwidth]{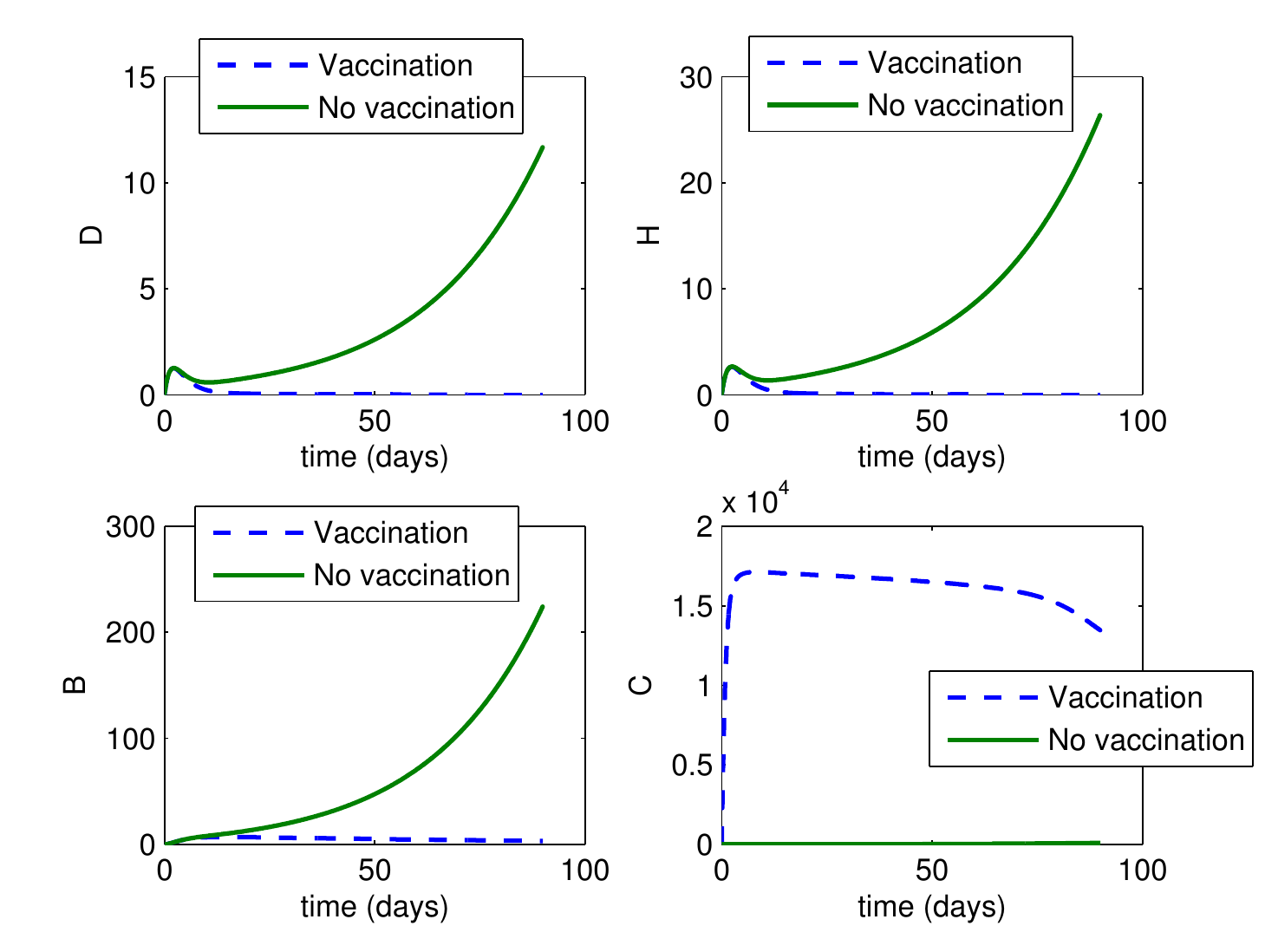}
\caption{Individuals $D(t)$, $H(t)$, $B(t)$ and $C(t)$, 
with the parameter values from 
Table~\ref{table:parameters}.
In dashed line, the case of vaccination without limit
on the supply of vaccines; in continuous line,
the case of no vaccination.}
\label{fig:DHBC:vac:novac}
\end{figure}

\begin{figure}[!htb]
\centering
\subfloat[\footnotesize{Optimal control $u(t)$}]{\label{control:Vac:NoLim}
\includegraphics[width=0.45\textwidth]{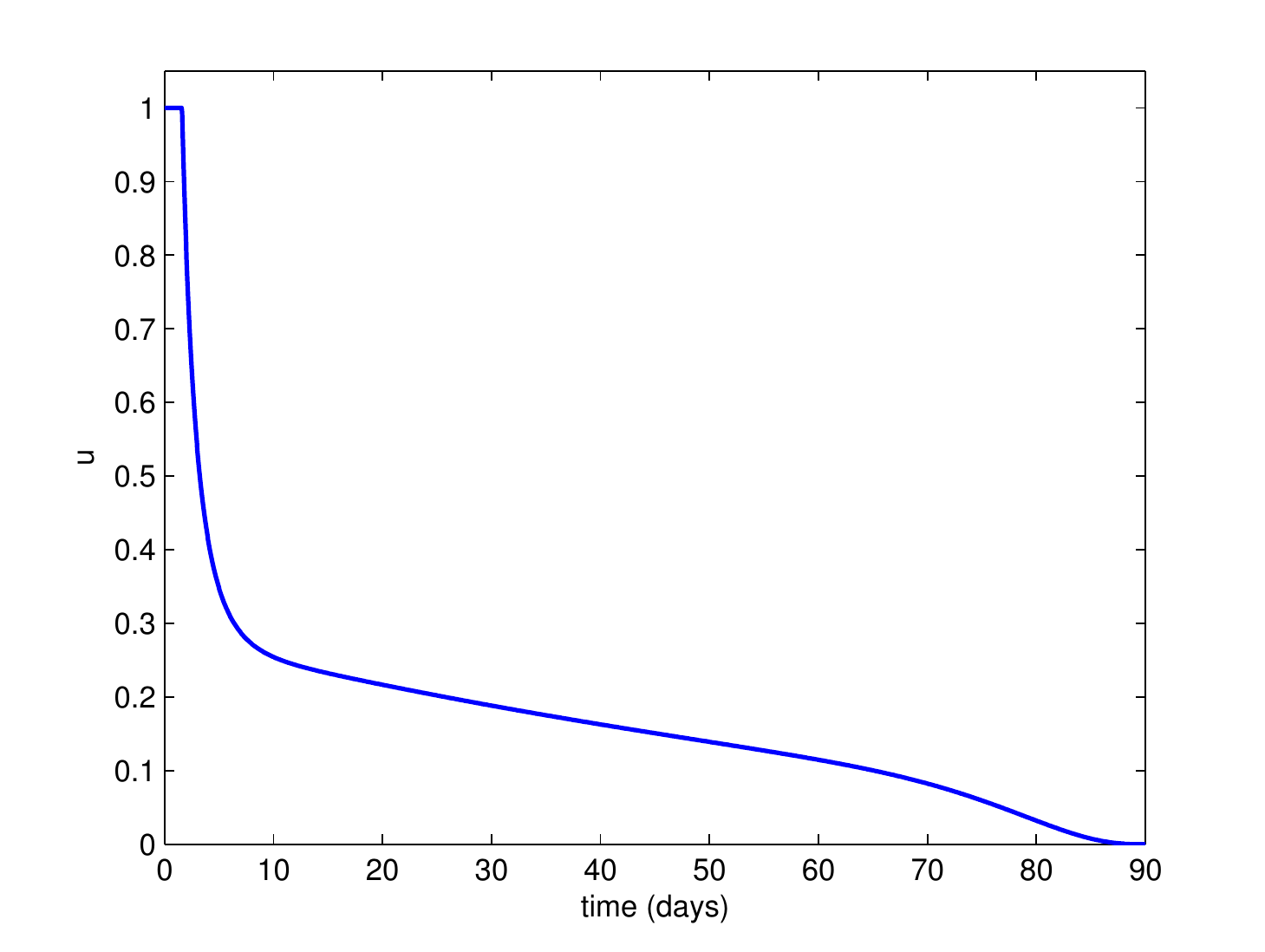}}
\subfloat[\footnotesize{Number of vaccines $W(t)$}]{\label{W:NoLim}
\includegraphics[width=0.45\textwidth]{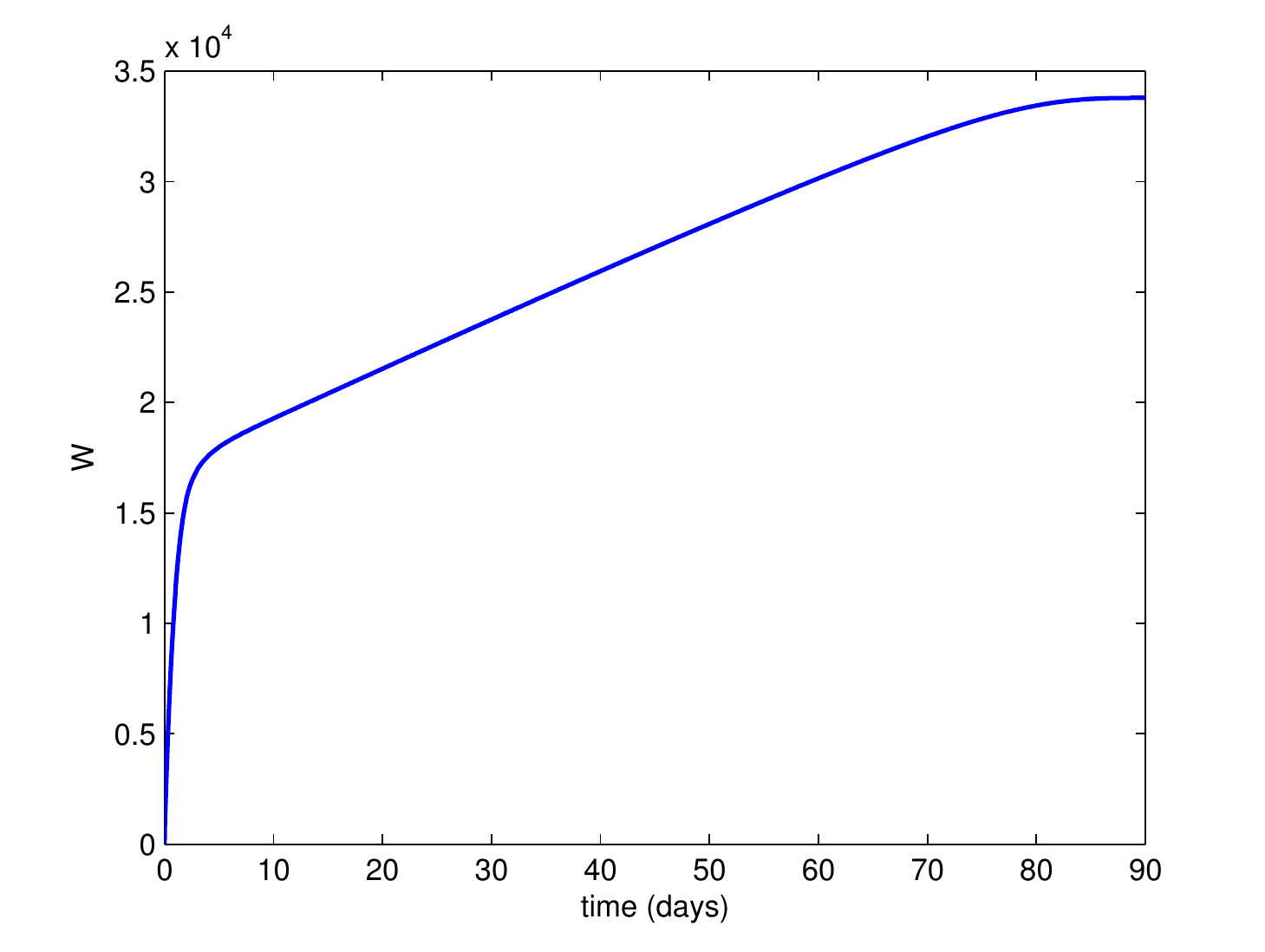}}
\caption{Optimal control and number of vaccines with 
the parameter values from Table~\ref{table:parameters}, 
when an unlimited supply of vaccines is available.}
\label{fig:control:vacc:NoLim}
\end{figure}

If we compare the case where there is no vaccination with the opposite case 
of unlimited supply of vaccines, we observe that at the end of 90 days the 
class of completely recovered individuals has approximately $86.5$ individuals 
in the case of no vaccination and $13468$ in the case of unlimited vaccination, 
which represents $74.82$ per cent of the total population 
(see Fig.~\ref{fig:SEIR:vac:novac}--\ref{fig:DHBC:vac:novac}). If vaccines are 
available, then the number of individuals that develop active disease is less 
than one at the end of 7.5 days and less that $0.1$ at the end of $32.4$ days. 
In the case of no vaccination, the class of active infected individuals 
has $61.7$ individuals at the end of 90 days (see Fig.~\ref{fig:SEIR:vac:novac}). 
If no vaccination is provided, then the number of deaths, hospitalizations 
and burials increases from $1.2$ to $262.6$, when compared to the case of 
unlimited supply of vaccines (see Fig.~\ref{fig:DHBC:vac:novac}).%
The optimal vaccination policy suggested by the solution of the optimal problem, 
implies a vaccination of 100 per cent of the susceptible population for 
approximately $1.62$ days followed by a fast reduction of the fraction of 
susceptible population that is vaccinated. This is based on the fact that 
the vaccine is effective and once all the susceptible population is vaccinated 
in a short period of time, then the number of susceptible individuals immediately 
decreases, since they are transferred to the class of completely recovered 
individuals, as well as the need of vaccination (see Fig.~\ref{control:Vac:NoLim}). 
The previous results show the importance of an effective vaccine for Ebola virus 
and the very good results that can be attained if the number of available 
vaccines satisfies the needs of the population.


\subsection{Limited total number of vaccines}
\label{sec:6.2}

\begin{figure}[!htb]
\centering
\subfloat[\footnotesize{Cumulative confirmed cases}]{\label{CumInfLimFinal}
\includegraphics[width=0.45\textwidth]{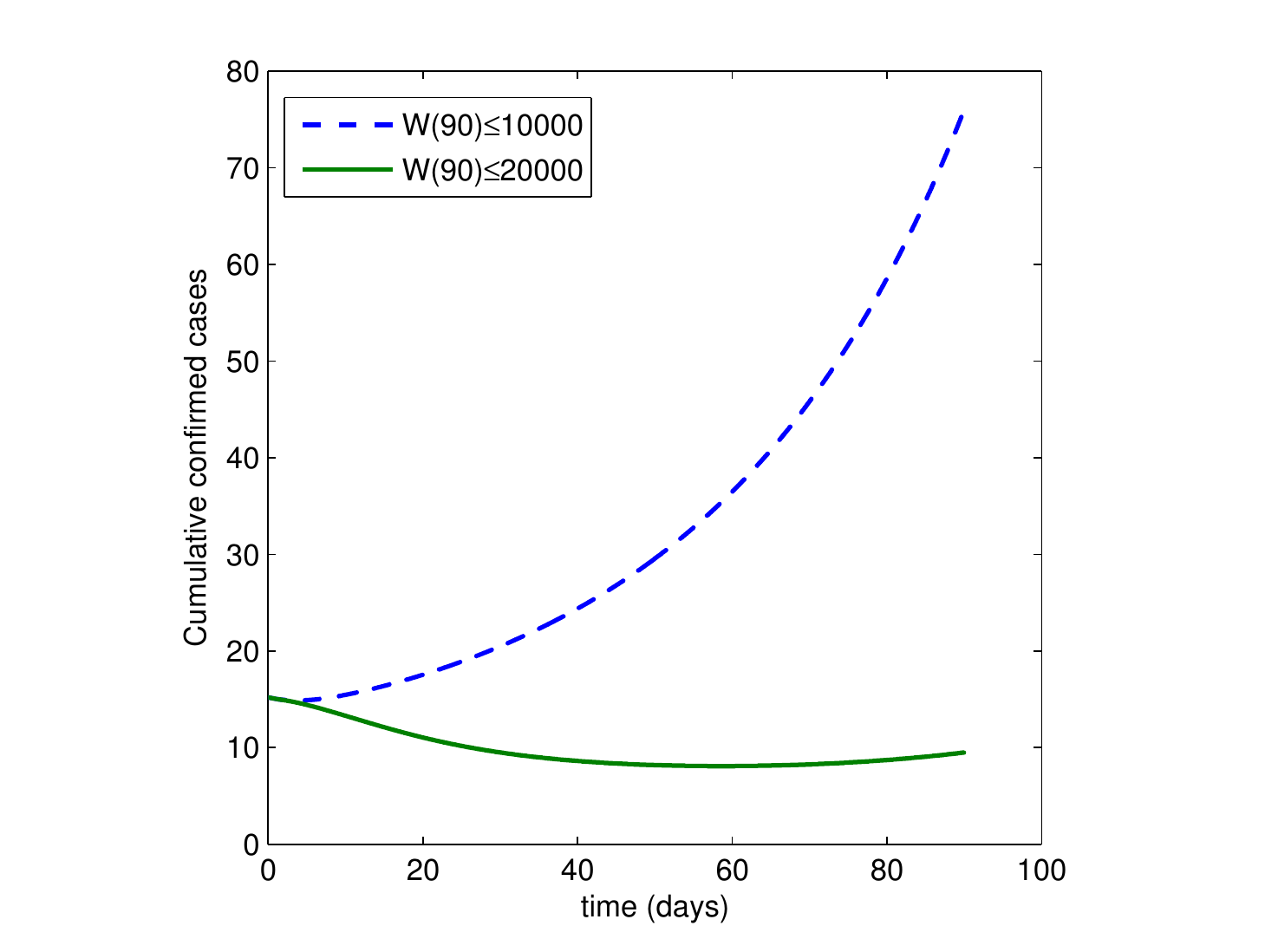}}
\subfloat[\footnotesize{Optimal control}]{\label{control:Vac:Lim10m20m}
\includegraphics[width=0.45\textwidth]{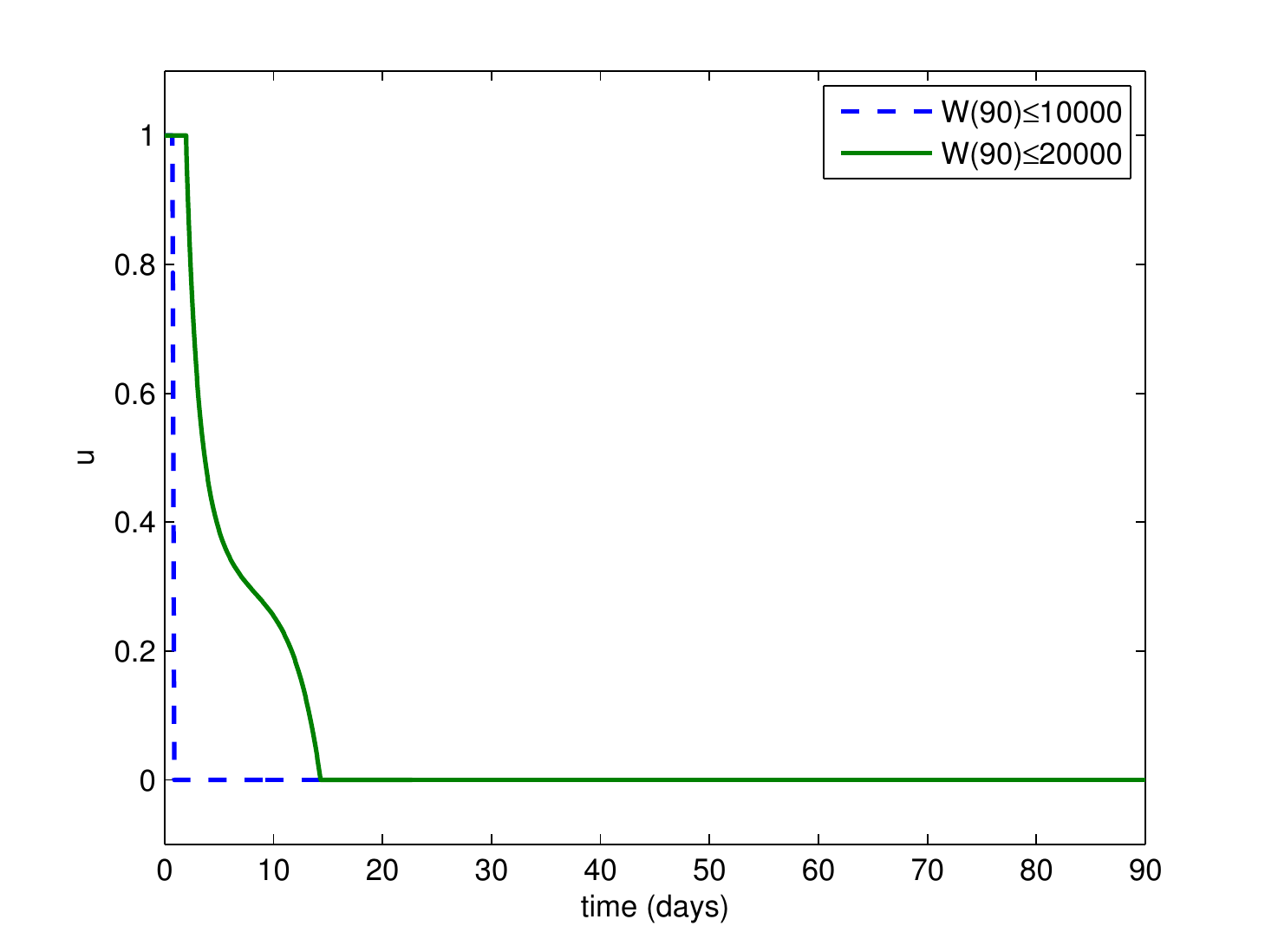}}
\caption{(a) Cumulative confirmed cases. (b) Optimal control for the case
of limited total number of vaccines. Dashed line for 
$W(90)\leq 10000$ and continuous line for $W(90)\leq 20000$.}
\label{fig:total:I:Vac:10mil20mil}
\end{figure}

In Fig.~\ref{W:NoLim}, we observe that at the end of 90 days, $33786$ 
vaccines were used, if the supply of vaccines is unlimited. In this section, 
we consider the case where the total number of vaccines used in the 90 days 
period is limited. We consider the case where the total number of vaccines 
available is lower or equal than the initial number of susceptible individuals 
($W(90)\leq 10000$, $W(90)\leq 11000$, $W(90)\leq 13000$, $W(90)\leq 15000$, 
$W(90)\leq 16000$, $W(90)\leq 18000$) and the case where the total number 
of vaccines available is bigger than the initial number of susceptible individuals 
($W(90) \leq 20000$). We first consider $w_1 = w_2 = 1$. The cumulative confirmed 
cases (see Fig.~\ref{fig:total:I:Vac:10mil20mil}) increases in time in the case 
$W(90) \leq 10000$ and decreases in the case $W(90) \leq 20000$, $t \in [0, 90]$. 
At the end of the 90 days period, the total number of individuals who got active 
infection is approximately $76$ and $9.5$ individuals in the case 
$W(90) \leq 10000$ and $W(90) \leq 20000$, respectively.
\begin{figure}[!htb]
\centering
\includegraphics[width=0.90\textwidth]{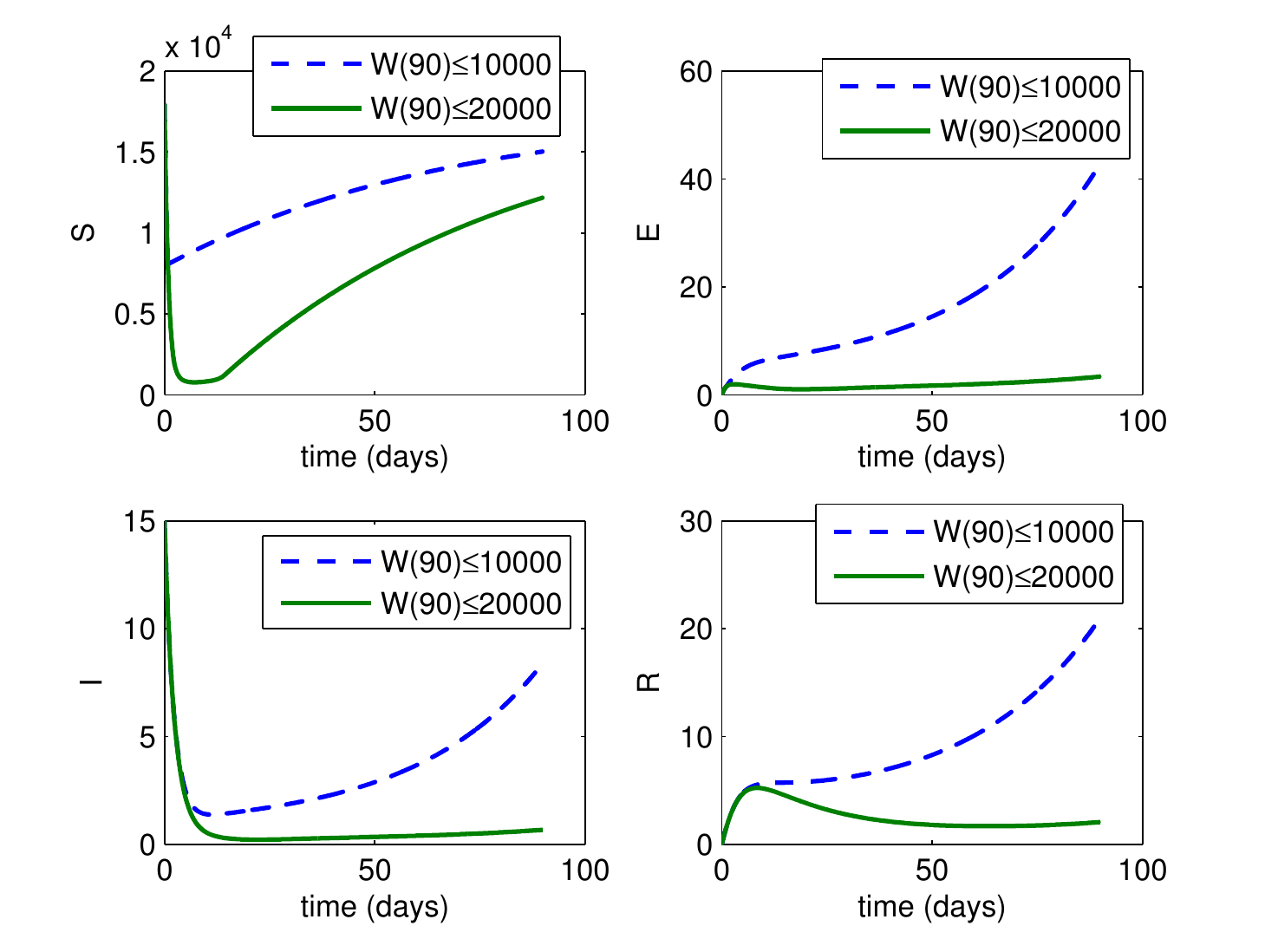}
\caption{Individuals $S(t)$, $E(t)$, $I(t)$ and $R(t)$. 
The dashed line represents the case where $W(90) \leq 10000$ 
and the continuous line represents the case where $W(90) \leq 20000$.}
\label{fig:SEIR:vac:Limit:10m20m}
\end{figure}
\begin{figure}[!htb]
\centering
\includegraphics[width=0.90\textwidth]{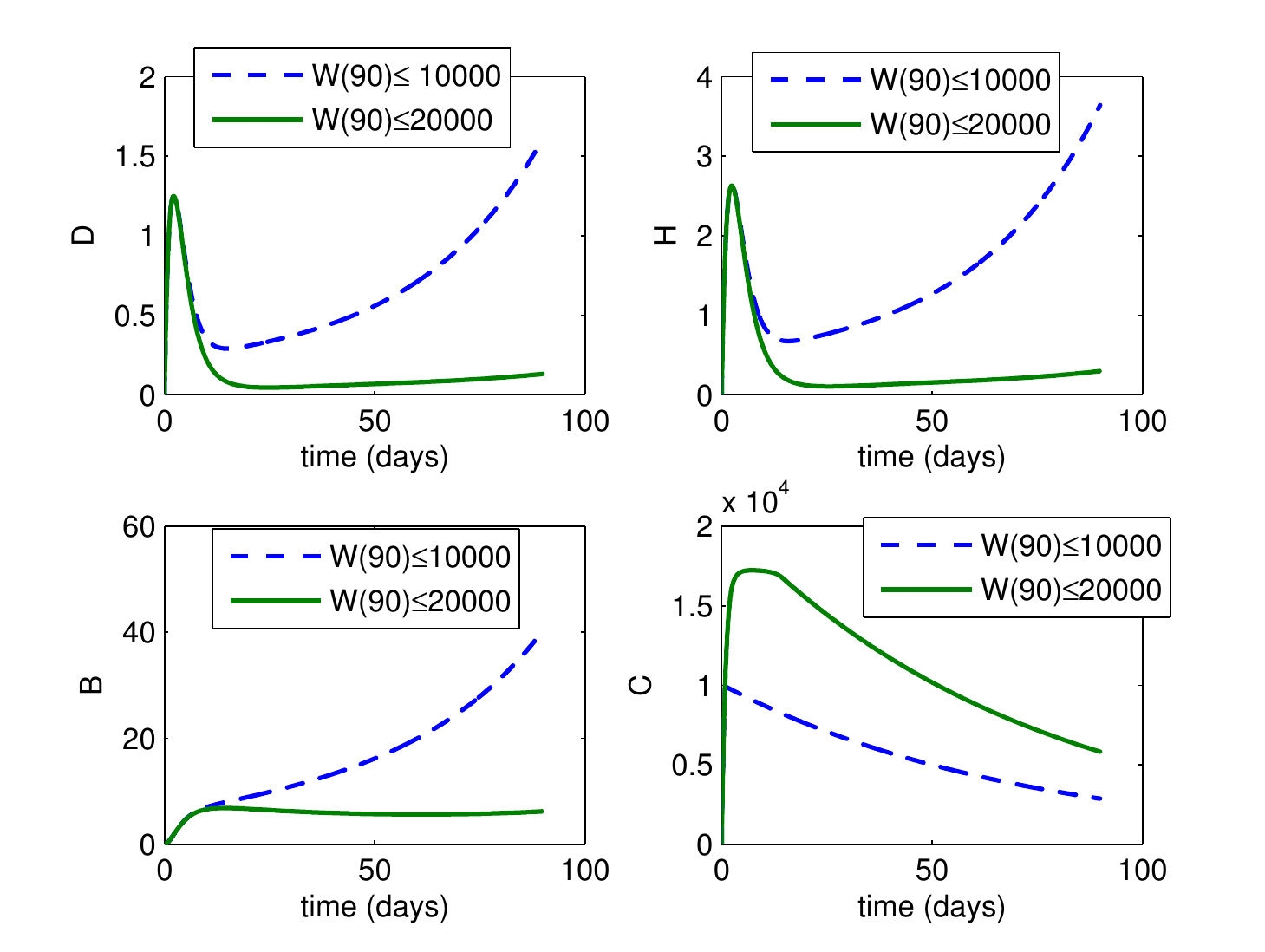}
\caption{Individuals $D(t)$, $H(t)$, $B(t)$ and $C(t)$. The dashed 
line represents the case where $W(90) \leq 10000$ and 
the continuous line represents the case where $W(90) \leq 20000$.}
\label{fig:DHBC:vac:Limit:10m20m}
\end{figure}
In the case $W(90)\leq 10000$, the optimal control takes the maximum value for 
less than one day (approximately 0.72 days) with a cost equal to $322.74$, and 
in the case $W(90) \leq 20000$, the optimal control takes the maximum value 
for approximately 2.2 days with a cost equal to $72.35$. The cost associated 
to the case $W(90) \leq 20000$ is lower than the one in the case $W(90) \leq 10000$, 
although more individuals are vaccinated, since the number of individuals in the 
class $I$ is lower. Namely, in the case $W(90)\leq 10000$, the number of individuals 
with active infection at the end of 90 days is equal to $I(90) = 8.4$ and in the 
case $W(90)\leq 20000$ the respective number is equal to $I(90)=0.67$. This means 
that in the case $W(90)\leq 10000$, in a epidemiological scenario corresponding 
to a basic reproduction number greater than one, 10000 vaccines will not be enough 
to eradicate the disease. Additionally, if we consider the maximum value for the 
total number of vaccines used during the period of 90 days to be equal to 
$11000$, $13000$, $15000$, $16000$ and $18000$, then we observe that the optimal 
control $u$ remains more time at the maximum value $1$ when the supply 
of vaccines is bigger, which means that when the total number of available 
vaccines is increased there will be resources to vaccinate all susceptible 
individuals for a longer period of time, which implies a bigger reduction 
of the number of individuals who get infected by the virus 
(see Fig.~\ref{control:Vac:variar} and \ref{control:Vac:variar:Zoom} for the 
optimal control strategy and respective zoom in the period of vaccination).
\begin{figure}[!htb]
\centering
\subfloat[\footnotesize{Optimal control}]{\label{control:Vac:variar}
\includegraphics[width=0.33\textwidth]{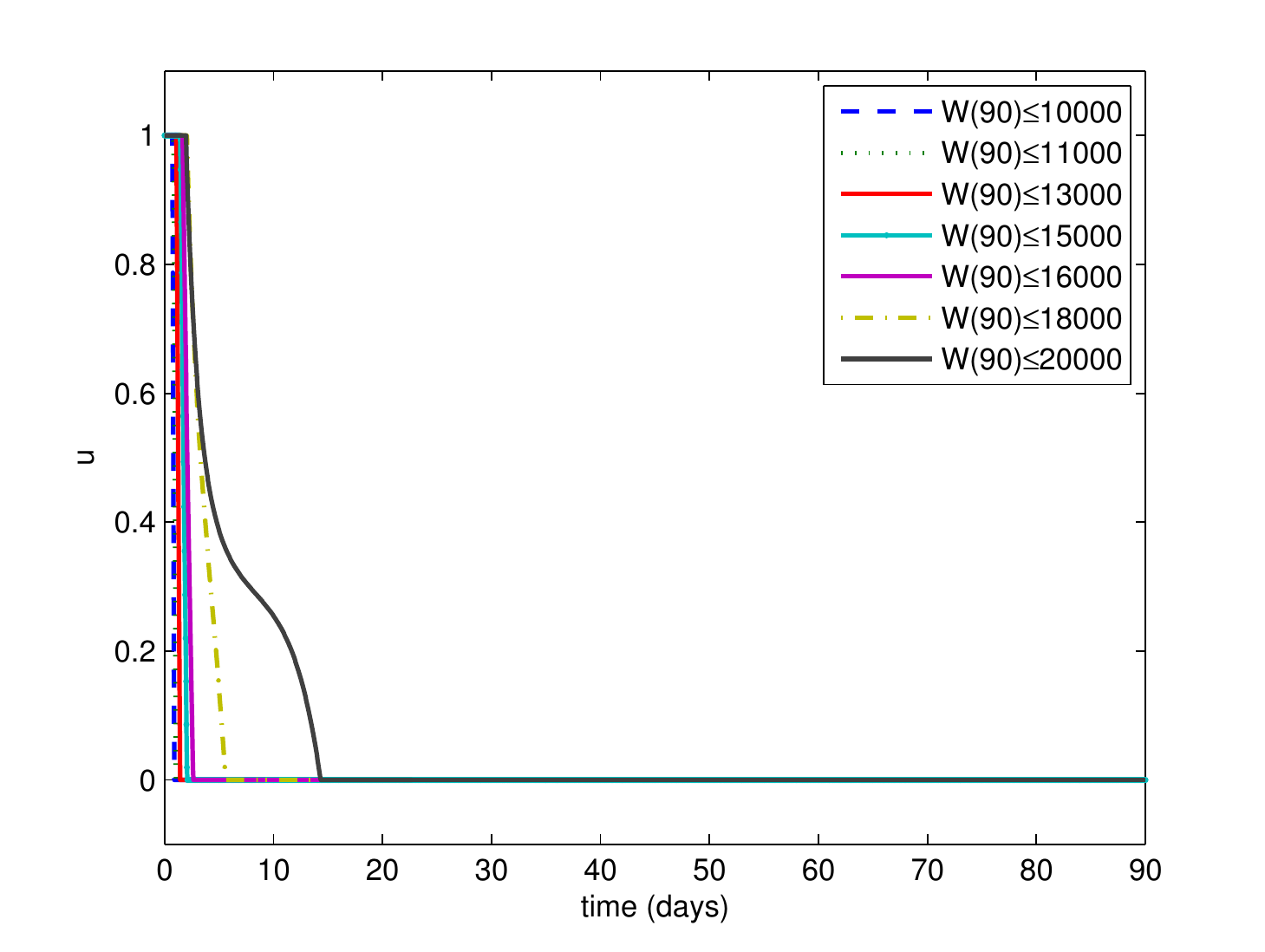}}
\subfloat[\footnotesize{Optimal control}]{\label{control:Vac:variar:Zoom}
\includegraphics[width=0.33\textwidth]{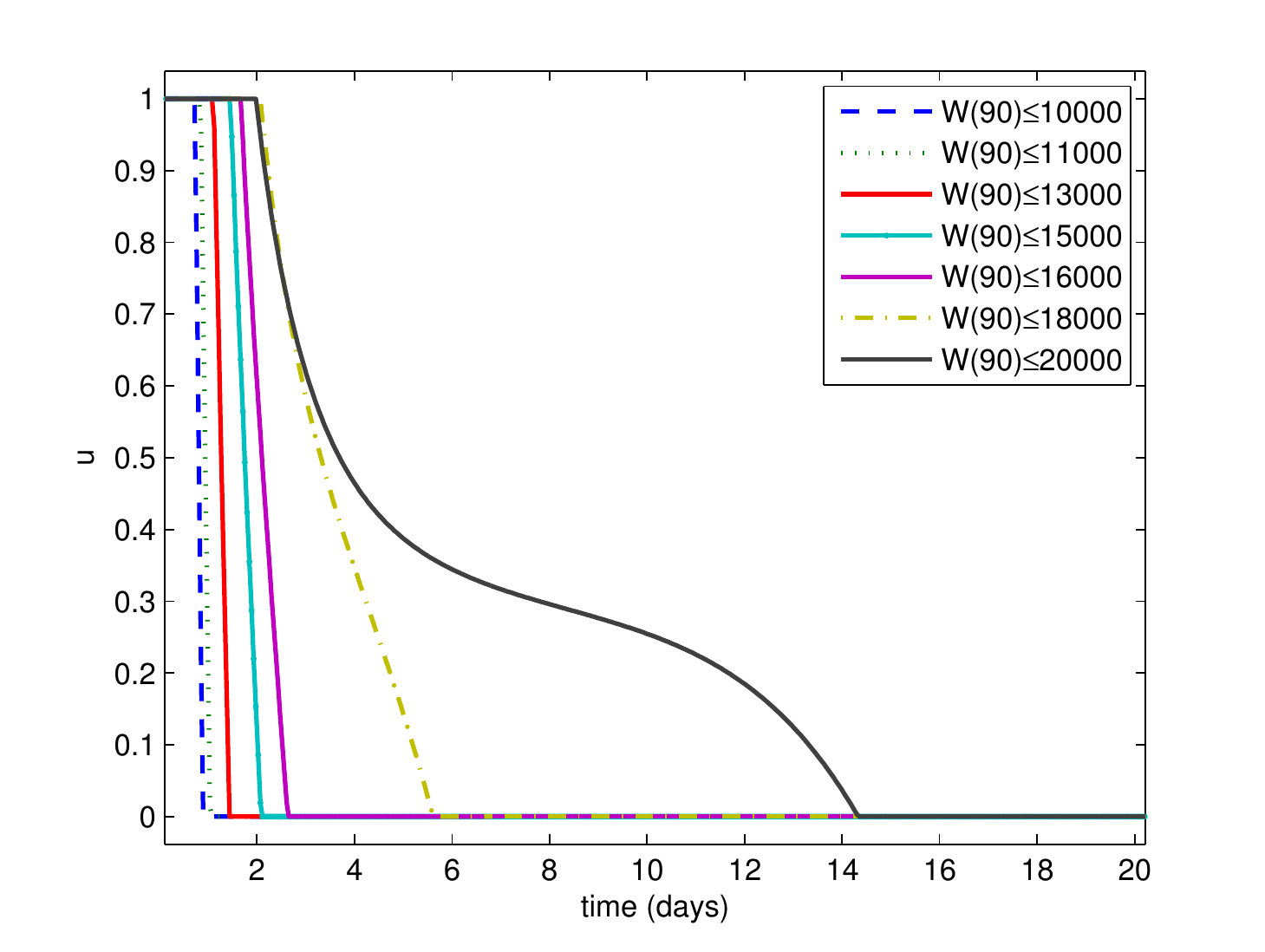}}
\subfloat[\footnotesize{Number of vaccines}]{\label{Vac:lim10m20m}
\includegraphics[width=0.33\textwidth]{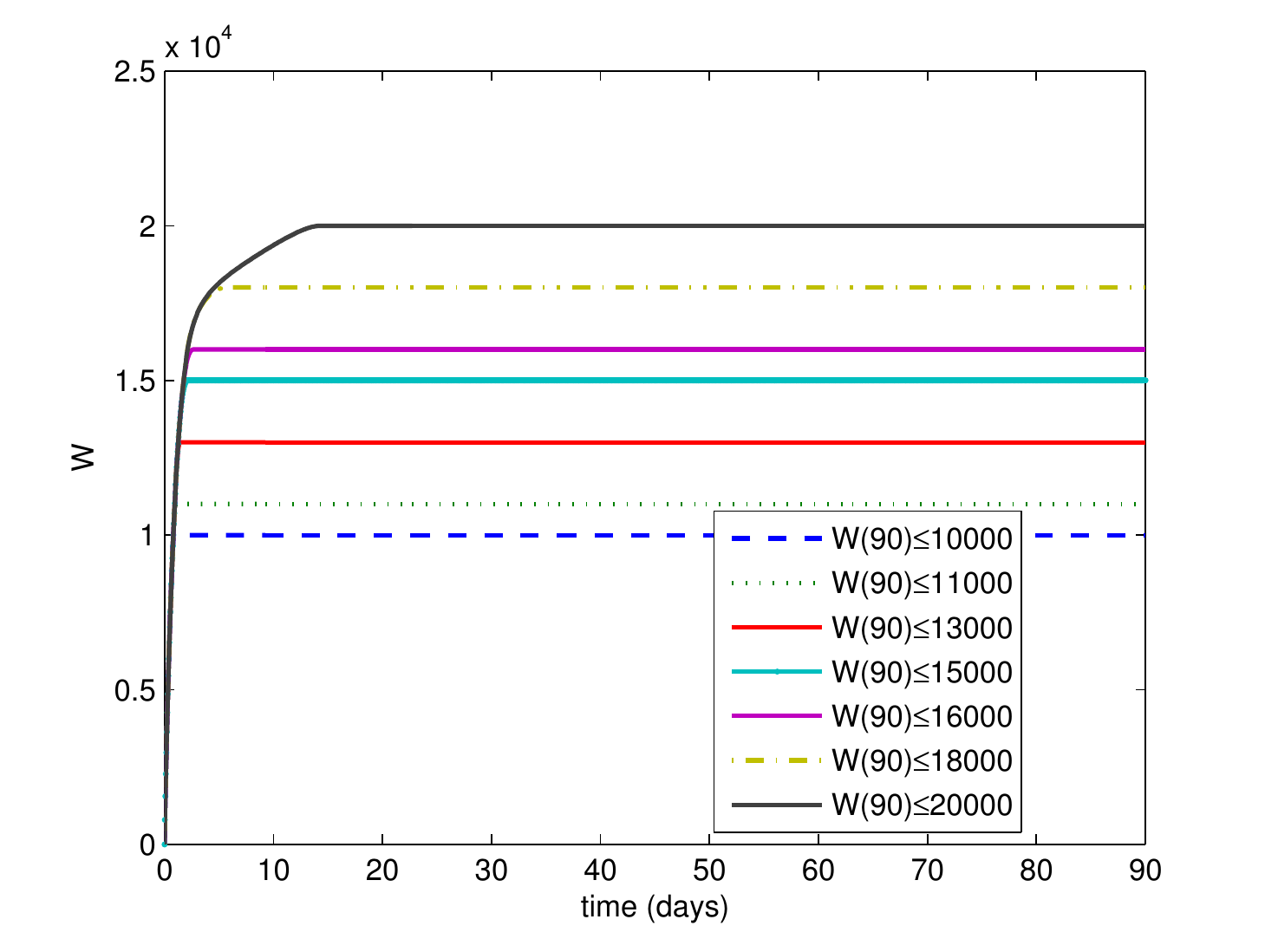}}
\caption{Optimal control $u(t)$ and number of vaccines $W(t)$
for $W(90)\leq 10000$, $W(90)\leq 11000$, $W(90)\leq 13000$,
$W(90)\leq 15000$, $W(90)\leq 16000$,
$W(90)\leq 18000$ and $W(90)\leq 20000$.}
\label{fig:control:vacc:limit:variar}
\end{figure}

\begin{figure}[!htb]
\centering
\subfloat[\footnotesize{Control}]{\label{control:Vac:LimA50A500}
\includegraphics[width=0.45\textwidth]{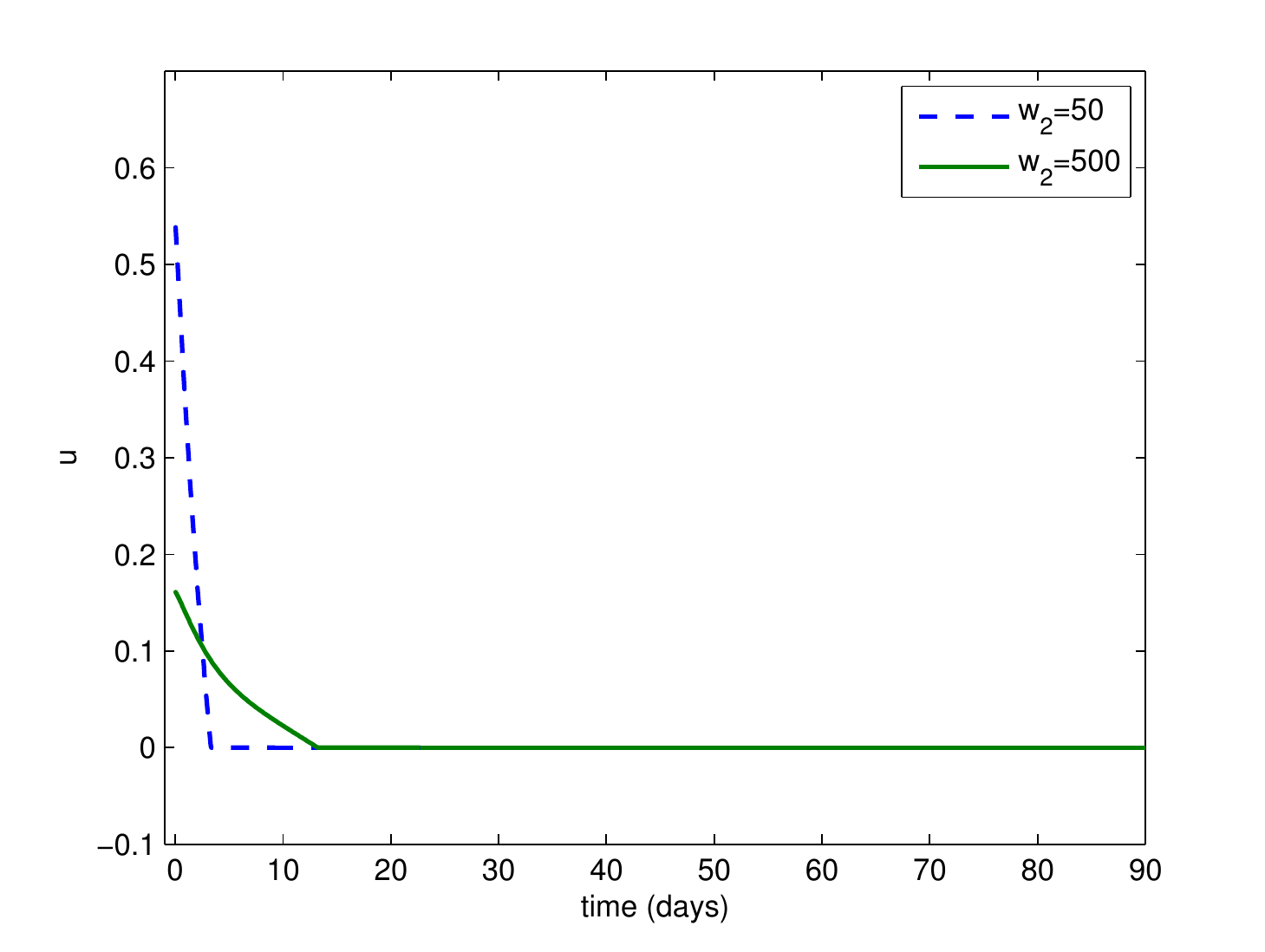}}
\subfloat[\footnotesize{Number of vaccines}]{\label{Vac:A50A500}
\includegraphics[width=0.45\textwidth]{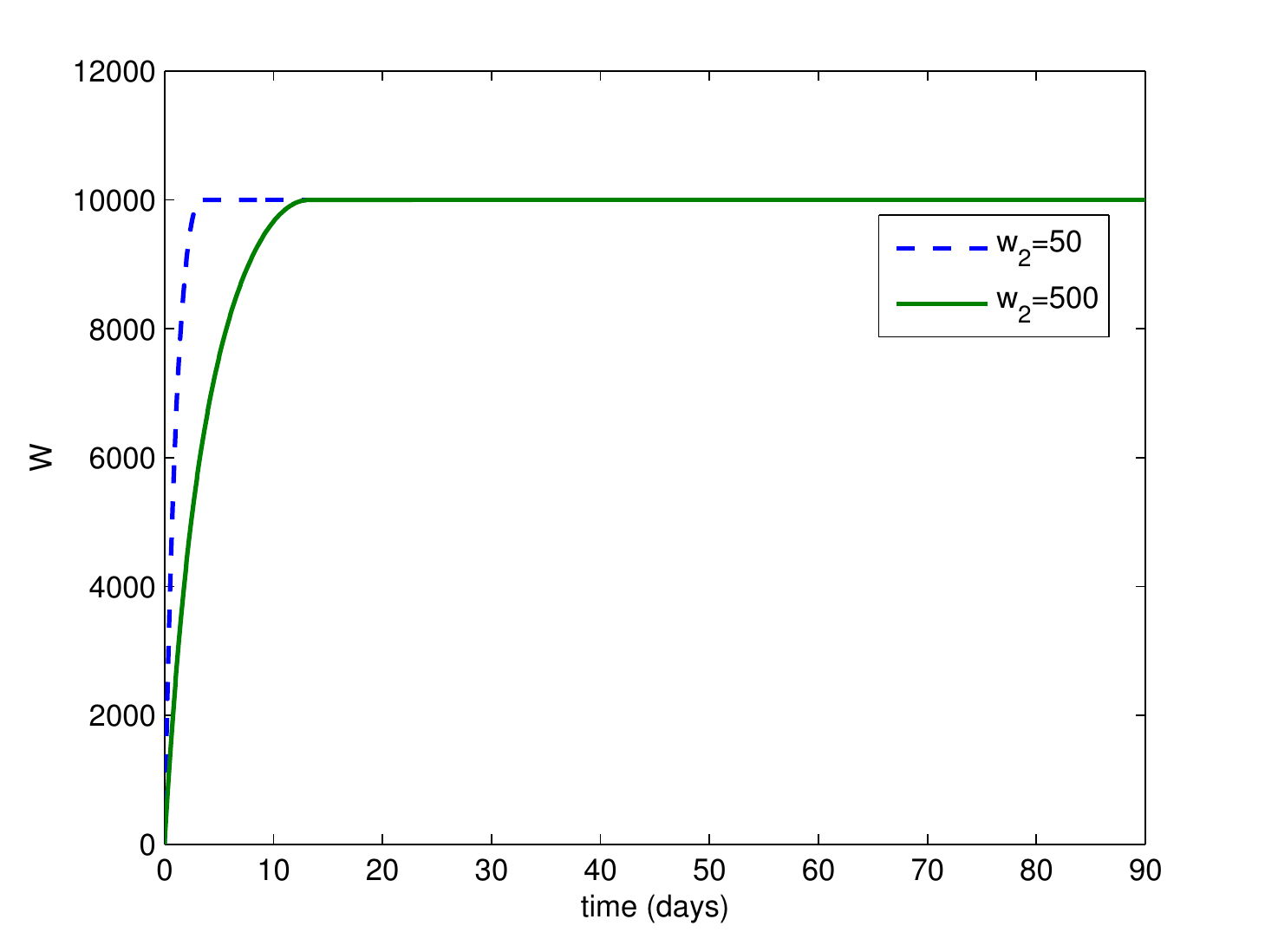}}
\caption{Optimal control $u(t)$ and number of vaccines $W(t)$
for $W(90)\leq 10000$. In dashed line the case $w_2=50$
and in continuous line the case $w_2 = 500$.}
\label{fig:control:vacc:limFinal:A50A500}
\end{figure}

Consider now the case where the weight constant associated with the cost of 
implementation of the vaccination strategy, designated by the optimal control 
$u$, is bigger than one, for example, consider $w_1 = 1$ and $w_2=50$, 
and $w_1 = 1$ and $w_2=500$. To simplify, consider in both cases $W(90)\leq 10000$ 
and $W(90)\leq 20000$. When we increase the weight constant $w_2$, 
the maximum value attained by the optimal control becomes lower than one 
(see Fig.~\ref{fig:control:vacc:limFinal:A50A500}). In the case 
$W(90)\leq 10000$ for $w_2 = 50$, the optimal control starts with the 
value $u(0) = 0.54$ and is a decreasing function with a cost function $344.3$. 
At the end of approximately 3.7 days, the control remains equal to zero. 
For $w_2 = 500$, the optimal control starts with the value $u(0) = 0.16$ 
and is also a decreasing function, with a cost $399.62$. At the end of 13.5 days, 
it remains equal to zero. The behavior of the optimal state variables
$S$, $E$, $I$, $R$, $D$, $H$, $B$ and $C$ are similar.


\subsection{Limited supply of vaccines at each instant of time}
\label{sec:6.3}

\begin{figure}[!htb]
\centering
\subfloat[\scriptsize{Total of active infected}]{\label{totalInfMixed}
\includegraphics[width=0.33\textwidth]{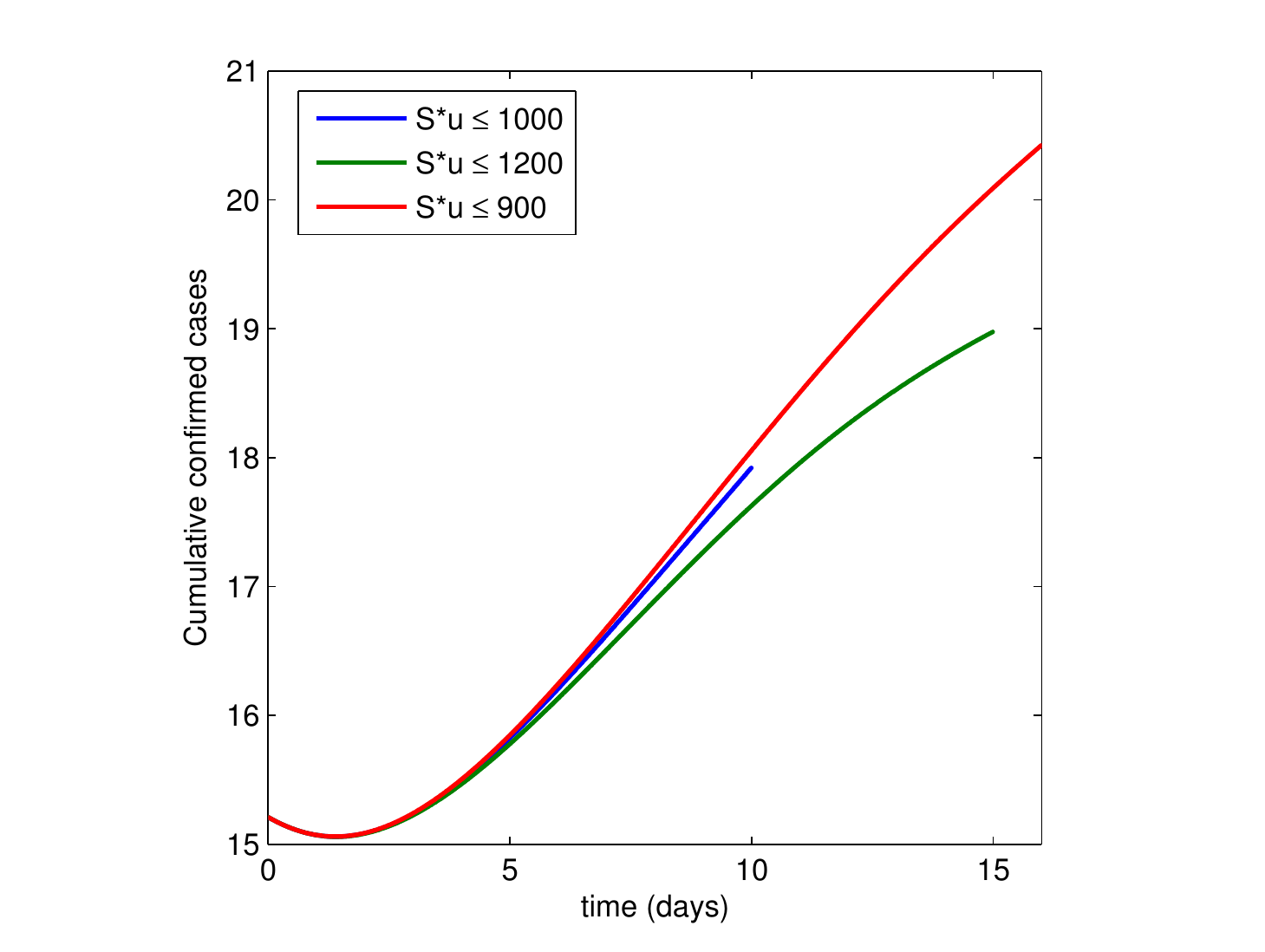}}
\subfloat[\scriptsize{Completely recovered}]{\label{Cmixed}
\includegraphics[width=0.34\textwidth]{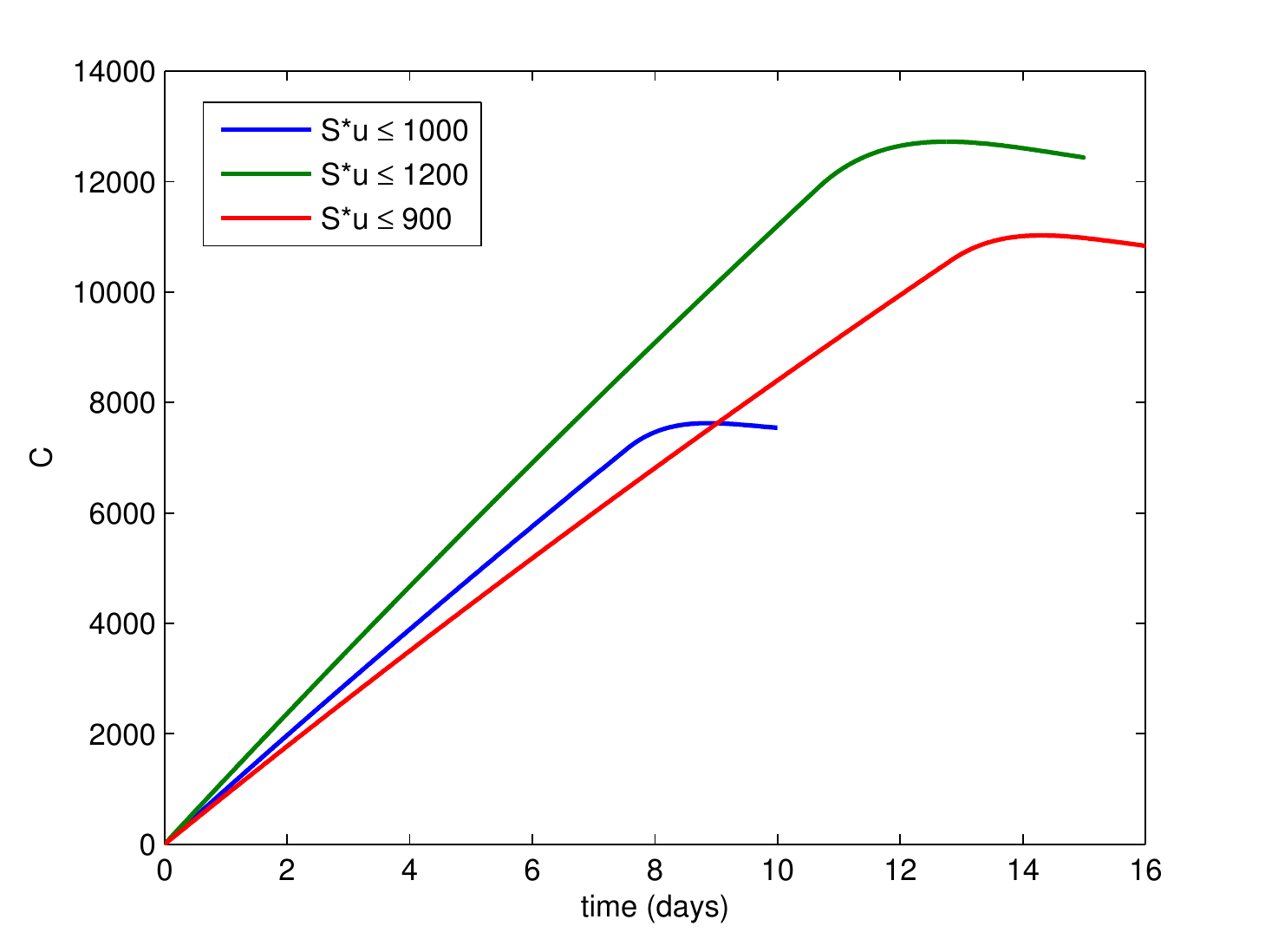}}
\subfloat[\scriptsize{Optimal control}]{\label{controlMixed}
\includegraphics[width=0.33\textwidth]{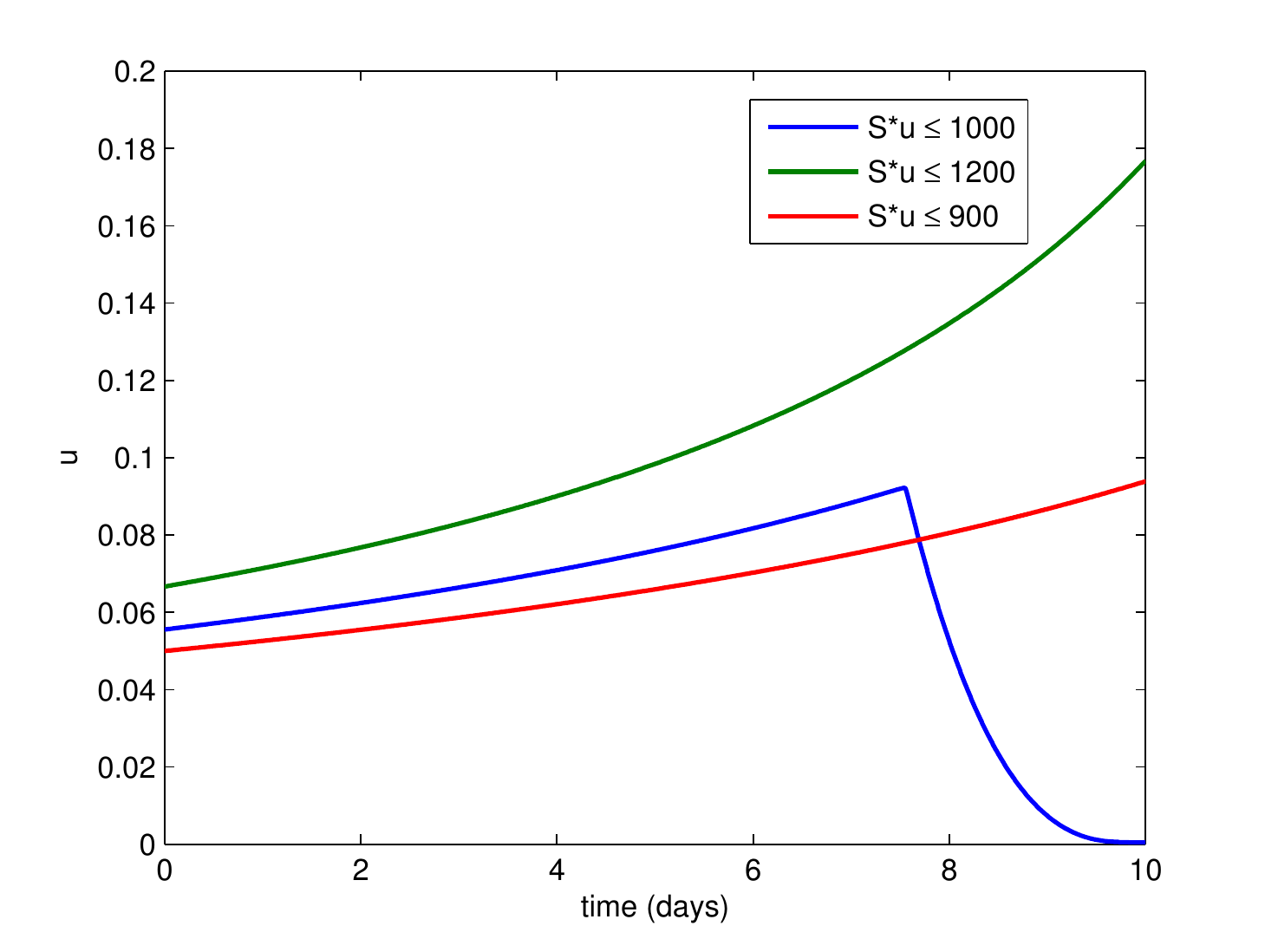}}
\caption{(a) Cumulative confirmed cases, (b) completely recovered,
(c) optimal control. In (a), (b) and (c) the following mixed constraints
are considered: $S(t) u(t) \leq 1000$ for all
$t \in [0, 10]$, $S(t) u(t) \leq 1200$ for all
$t \in [0, 15]$, and $S(t) u(t) \leq 900$ for all $t \in [0, 16]$.}
\label{fig:total:inf:mixed}
\end{figure}

From previous results, we observe that when there is no limit on the supply 
of vaccines, or when the total number of used vaccines is limited, the optimal 
vaccination strategy implies a vaccination of 100 per cent of the susceptible 
population in a very short period of time, sometimes smaller than one day. 
But we know that in practice this is a very difficult task, since there are 
limitations in the number of vaccines available and also in the number of health 
care workers or humanitarian teams in the regions affected by Ebola virus 
with capacity to vaccinate such a big number of individuals almost simultaneously. 
{}From this point of view, it is important to study the case where there is a 
limited supply of vaccines at each instant of time. In this section, we consider 
$w_1 = w_2 = 1$, a shorter interval of time $[0, t_f]$, with $t_f = 10$, $15$, 
$16$, and we assume that at each instant of time there exist only 1000, 1200 
and 900 available vaccines, respectively. From our point of view, these numbers 
of available vaccines at each instant of time and the number of days considered, 
correspond to possible real scenarios, which are possible to implement 
in a concrete endemic region and at the same time characterize lack of human 
and material resources to vaccinate the susceptible population in a short period 
of time. From the numerical simulations, for such mixed constraints, the number 
of cumulative confirmed cases increases with time (see Figure~\ref{totalInfMixed}). 
The cost associated with the vaccination campaign, associated with the solution 
of the optimal control problem with the mixed constraint $S(t) u(t) \leq 1000$, 
$t \in [0, 10]$, is equal to $45.8879$. Such solution is the less costly 
of the three considered, followed by the constraint $S(t) u(t) \leq 1200$ 
for $t \in [0, 15]$ with a cost of $55.079$. The most expensive vaccination 
strategy is the one associated with the mixed constraint $S(t) u(t) \leq 900$, 
$t \in [0, 16]$, with a cost of $59.109$. The strategy associated with the 
constraint $S(t) u(t) \leq 1000$ is the one where a lowest number of susceptible 
individuals completely recover through vaccination, with $7540.9$ individuals 
in the class $C$ at the end of 10 days. If we consider that at each instant 
of time there are 1200 vaccines available during a period of 15 days, 
then $12438$ completely recover. This is the strategy with more individuals 
in the class $C$. If we consider 16 days, but only 900 vaccines available 
for each instant of time, then only $10839$ individuals completely recover 
(see Fig.~\ref{Cmixed}). For all three mixed constraint situations, the number 
of individuals in the classes $E$, $I$, $R$, $D$, $H$ and $B$ does not change 
significantly (therefore, the figures with these classes are omitted). As the number 
of available vaccines represent a small percentage of the susceptible population, 
in the three cases the optimal vaccination strategies for the constraints 
$S(t) u(t) \leq 1200$ and $S(t) u(t) \leq 900$ suggest that the percentage 
of the susceptible population that is vaccinated is always inferior than 
18 percent. In the case of the constraint $S(t) u(t) \leq 1000$, this percentage 
is always inferior to 8 percent (see Fig.~\ref{controlMixed}).


\section{Discussion}
\label{sec:7}

We assume that, in a near future, an effective vaccine against 
the Ebola virus will be available. Under this assumption, 
three different scenarios have been studied: unlimited supply of vaccines; 
limited total number of vaccines to be used; and limited supply 
of vaccines at each instant of time. We have solved the optimal 
control problems analytically and we have performed a number 
of numerical simulations in the three aforementioned vaccination scenarios.

Some authors have already considered the optimal control problem 
with vaccination for Ebola disease, but always with unlimited 
supply of vaccines \cite{Dure,MR3349757}. It turns out that 
the solution to this mathematical problem is obvious: the solution 
consists to vaccinate all susceptible individuals in the beginning 
of the outbreak. This is a very particular case of our work, 
investigated in Section~\ref{sec:6.1} (see Figure~\ref{fig:control:vacc:NoLim}). 
If vaccines are available without any restriction, then one could  
completely eradicate Ebola in a very short period of time. These results 
show the importance of an effective vaccine for Ebola virus 
and the very good results that can be attained if the number 
of available vaccines satisfy the needs of the population. 
Unfortunately, such situation is not realistic: in case  
an effective vaccine for Ebola virus will appear, there always 
will be restrictions on the number of available vaccines  
as well as constraints on how to inoculate them in a proper way 
and in a short period of time; economic problems might also exist.

In our work, for first time in the literature of Ebola,
an optimal control problem with state and control 
constraints has been considered. Mathematically, it represents a health public 
problem of limited total number of vaccines. The results obtained 
in Section~\ref{sec:6.2} provide useful information on the number 
of vaccines to be bought, in order to reduce the number 
of new infections with minimum cost. For example, the results between 
10000 and 20000 vaccines (in 90 days) are completely different. 
With 10000 vaccines, the number of cumulative infected cases continues 
to increase, while with 20000 vaccines it is already possible 
to decrease the new infections. The optimal solution, in this case, 
is similar to the case of unlimited supply of vaccines, that is, 
it implies a vaccination of 100 per cent of the susceptible population 
in a very short period of time. In practice, this is an unrealistic task, 
due to the necessary number of vaccines and humanitarian teams 
in the regions affected by Ebola. Therefore, we conclude that 
it is important to study the case where there is a limited supply 
of vaccines at each instant of time. This was investigated 
in Section~\ref{sec:6.3}. This situation is much richer 
and the optimal control solution is not obvious. For a given number 
of available vaccines at each instant of time, we have a different solution, 
which is the optimal rate of susceptible individuals that should be vaccinated. 
In this case, the optimal control implies the vaccination of a small subset 
of the susceptible population. It remains the ethical problem 
of how to choose the individuals to be vaccinated.


\section*{Acknowledgments}

Three reviewers deserve special thanks for helpful and constructive comments.
The work of Area was partially supported by the Ministerio de Econom\'{\i}a 
y Competi\-tividad of Spain, under grants MTM2012--38794--C02--01 and MTM2016--75140--P, 
co-financed by the European Community fund FEDER. Nda\"{\i}rou acknowledges the 
AIMS-Cameroon 2014--2015 fellowship, Nieto the partial financial 
support by the Ministerio de Econom\'{\i}a y Competitividad of Spain under 
grants MTM2010--15314, MTM2013--43014--P and MTM2016--75140--P, and Xunta 
de Galicia under grants R2014/002 and GRC 2015/004, co-financed by the 
European Community fund FEDER. Silva was supported through the Portuguese 
Foundation for Science and Technology (FCT) post-doc fellowship SFRH/BPD/72061/2010.
The work of Silva and Torres was partially supported by FCT through CIDMA 
and project UID/MAT/04106/2013, and by project TOCCATA, 
reference PTDC/EEI-AUT/2933/2014, funded by Project 
3599 -- Promover a Produ\c{c}\~ao Cient\'{\i}fica e Desenvolvimento
Tecnol\'ogico e a Constitui\c{c}\~ao de Redes Tem\'aticas (3599-PPCDT)
and FEDER funds through COMPETE 2020, Programa Operacional
Competitividade e Internacionaliza\c{c}\~ao (POCI), and by national
funds through FCT.



\medskip
Received February 2016; revised November 2016; accepted March 2017.
\medskip



\begin{thebibliography}{10}

\bibitem{Dure} [10.1186/s40249-016-0161-6]
\newblock M. D. Ahmad, M. Usman, A. Khan and M. Imran,
\newblock Optimal control analysis of Ebola disease 
with control strategies of quarantine and vaccination,
\newblock {\emph{Infect. Dis. Poverty}} \textbf{5} (2016), no.~72, 12~pp. 

\bibitem{Althaus} [10.1371/currents.outbreaks.91afb5e0f279e7f29e7056095255b288]
\newblock C. L. Althaus, 
\newblock Estimating the reproduction number of Ebola 
virus (EBOV) during the 2014 outbreak in west Africa,
\newblock {\emph{PLOS Currents Outbreaks}}, Edition 1, 2014. 

\bibitem{MR3394468} (MR3394468) [10.1186/s13662-015-0613-5]
\newblock I. Area, H. Batarfi, J. Losada, J. J. Nieto, W. Shammakh and A. Torres,
\newblock On a fractional order {E}bola epidemic model,
\newblock \emph{Adv. Difference Equ.} \textbf{2015} (2015), Art. ID 278, 12~pp.

\bibitem{Area:in:press} [10.1002/mma.3794]
\newblock I. Area, J. Losada, F. Nda\"{\i}rou, J. J. Nieto and D. D. Tcheutia,
\newblock Mathematical modeling of 2014 Ebola outbreak,
\newblock \emph{Math. Method. Appl. Sci.}, in press.

\bibitem{atangana2014} [10.1155/2014/261383]
\newblock A. Atangana and E. F. Doungmo~Goufo,
\newblock On the mathematical analysis of Ebola hemorrhagic fever: 
deathly infection disease in west african countries,
\newblock \emph{BioMed Research International} 
\textbf{2014} (2014), Art. ID 261383, 7~pp.

\bibitem{SEIR:Rosario:2014} (MR3181992) [10.3934/mbe.2014.11.761]
\newblock M. H. A. Biswas, L. T. Paiva and M. R. de~Pinho,
\newblock A {SEIR} model for control of infectious diseases with constraints,
\newblock \emph{Math. Biosci. Eng.} \textbf{11} (2014), no.~4, 761--784.

\bibitem{bwaka} [10.1086/514308]
\newblock M. A. Bwaka et al., 
\newblock Ebola hemorrhagic fever in Kikwit, Democratic Republic 
of the Congo: clinical observations in 103 patients, 
\newblock {\emph{J. Infect. Dis.}} \textbf{179} (1999), Suppl.~1, S1--S7.

\bibitem{Cesari_1983} (MR688142) [10.1007/978-1-4613-8165-5]
\newblock L. Cesari, 
\newblock \emph{Optimization---theory and applications},
Vol.~17 of Applications of Mathematics (New York),
\newblock Springer-Verlag, New York, 1983.

\bibitem{Chowell} (MR2082775) [10.1016/j.jtbi.2004.03.006]
\newblock G. Chowell, N. W. Hengartner, C. Castillo-Chavez, P. W. Fenimore and J. M. Hyman, 
\newblock The basic reproductive number of Ebola and the effects 
of public health measures: the cases of Congo and Uganda,
\newblock {\emph{J. Theor. Biol.}} \textbf{229} (2004), no.~1, 119--126.

\bibitem {ChowellN} [10.1186/s12916-014-0196-0]
\newblock G. Chowell and H. Nishiura,
\newblock Transmission dynamics and control of Ebola virus disease (EVD): a review,
\newblock {\emph{BMC Med.}} \textbf{12} (2014), no.~196, 16~pp. 

\bibitem{livro:Clarke:2013} (MR3026831) [10.1007/978-1-4471-4820-3]
\newblock F. Clarke, 
\newblock \emph{Functional analysis, calculus of variations and optimal control}, 
Vol. 264 of Graduate Texts in Mathematics,
\newblock Springer, London, 2013.

\bibitem{SIAM:Clarke:Rosario} (MR2683896) [10.1137/090757642]
\newblock F. Clarke and M. R. de~Pinho,
\newblock Optimal control problems with mixed constraints,
\newblock \emph{SIAM J. Control Optim.} \textbf{48} (2010), no.~7, 4500--4524.

\bibitem{diekmann} (MR1057044) [10.1007/BF00178324]
\newblock O. Diekmann, J. A. P. Heesterbeek and J. A. J. Metz, 
\newblock On the definition and the computation of the basic reproduction 
ratio ${R}_{0}$ in models for infectious diseases in heterogeneous populations,
\newblock \emph{J. Math. Biol.} \textbf{28} (1990), no.~4, 365--382.

\bibitem{dowell} [10.1086/514284]
\newblock S. F. Dowell,
\newblock Transmission of Ebola hemorrhagic fever: a study of risk factors 
in family members, Kikwit, Democratic Republic of the Congo, 1995, 
\newblock {\emph{J. Infect. Dis.}} {\textbf{179}} (1999), Suppl.~1, S87--S91.

\bibitem{OC:HIV:PLoSCompBio:2015} [10.1371/journal.pcbi.1004200]
\newblock S. Duwal, S. Winkelmann, C. Sch\"{u}tte and M. von Kleist,
\newblock Optimal treatment strategies in the context of 
'Treatment for Prevention' against HIV-1 in resource-poor settings,
\newblock \emph{PLoS Comput. Biol.} \textbf{11} (2015), no.~4, Art. ID e1004200, 30~pp.

\bibitem{Fleming_Rishel_1975} (MR0454768)
\newblock W. H. Fleming and R. W. Rishel,
\newblock \emph{Deterministic and stochastic optimal control},
\newblock Springer-Verlag, Berlin-New York, 1975.

\bibitem{Gire20141369} [10.1126/science.1259657]
\newblock S. K. Gire, A. Goba, K. G. Andersen, 
R. S. G. Sealfon, D. J. Park, L. Kanneh, et~al.,
\newblock Genomic surveillance elucidates Ebola virus origin 
and transmission during the 2014 outbreak,
\newblock \emph{Science}, \textbf{345} (2014), no.~6202, 1369--1372.

\bibitem{Hayden2014}
\newblock E. C. Hayden,
\newblock World struggles to stop Ebola,
\newblock \emph{Nature} \textbf{512} (2014), no.~7515, 355--356.

\bibitem{MyID:335} (MR3578107) [10.1080/15502287.2016.1231236] 
\newblock D. Hincapie-Palacio, J. Ospina and D. F. M. Torres,
\newblock Approximated analytical solution to an Ebola optimal control problem,
\newblock \emph{Int. J. Comput. Methods Eng. Sci. Mech.} \textbf{17} (2016), no.~5-6, 382--390. 
\newblock {\tt arXiv:1512.02843}

\bibitem{ebola:models}
\newblock J. Kaufman, S. Bianco and A. Jones,  
{\url{https://wiki.eclipse.org/Ebola_Models}}.

\bibitem{Kaushik2016254} [10.1016/j.bios.2015.08.040]
\newblock A. Kaushik, S. Tiwari, R. D. Jayant, A. Marty and M. Nair, 
\newblock Towards detection and diagnosis of Ebola virus disease at point-of-care,
\newblock \emph{Biosensors and Bioelectronics} \textbf{75} (2016), 254--272.

\bibitem{khan} [10.1086/514306]
\newblock A. S. Khan et al., 
\newblock The reemergence of Ebola hemorrhagic fever, Democratic Republic of the Congo, 1995,
\newblock {\emph{J. Infect. Dis.}} {\textbf{179}} (1999), Suppl.~1, S76--S86.

\bibitem{Ledzewicz:cancer:SIAM:2007} (MR2338438) [10.1137/060665294]
\newblock U. Ledzewicz and H. Sch\"{a}ttler, 
\newblock Antiangiogenic therapy in cancer treatment as an optimal control problem, 
\newblock \emph{SIAM J. Control Optim.} \textbf{46} (2007), no.~3, 1052--1079.

\bibitem{HYG:996392} [10.1017/S0950268806007217]
\newblock J. Legrand, R. F. Grais, P. Y. Boelle, A. J. Valleron and A. Flahault,
\newblock Understanding the dynamics of {E}bola epidemics,
\newblock \emph{Epidemiol. Infect.} \textbf{135} (2007), no.~4, 610--621.

\bibitem{lekone} [10.1111/j.1541-0420.2006.00609.x]
\newblock P. E. Lekone and B. F. Finkenst\"{a}dt,
\newblock Statistical inference in a stochastic epidemic {SEIR} model 
with control intervention: Ebola as a case study,
\newblock \emph{Biometrics} \textbf{62} (2006), no.~4, 1170--1177.

\bibitem{OC:HepatiticC:PLoSOne:2011} [10.1371/journal.pone.0022309]
\newblock N. K. Martin, A. B. Pitcher, P. Vickerman, A. Vassall and M. Hickman, 
\newblock Optimal control of hepatitis C antiviral treatment programme 
delivery for prevention amongst a population of injecting drug users,
\newblock \emph{PLoS One} \textbf{6} (2011), no.~8, e22309, 17~pp.

\bibitem{NeilanLenhart2010} (MR2744727) 
\newblock R. Miller~Neilan and S. Lenhart, 
\newblock \emph{An introduction to optimal control with an application in disease modeling}.
\newblock In: Modeling paradigms and analysis of disease transmission models,
Vol.~75 of DIMACS Ser. Discrete Math. Theoret. Comput. Sci. Amer. Math. Soc.,
Providence, RI, 2010, 67--81.

\bibitem{ndambi} [10.1086/514297]
\newblock R. Ndambi et al., 
\newblock Epidemiologic and clinical aspects of the Ebola virus epidemic in Mosango,
Democratic Republic of the Congo, 1995, 
\newblock {\emph{J. Infect. Dis.}} {\textbf{179}} (1999), Suppl.~1, S8--S10.

\bibitem{Ngwa} (MR3538903) [10.1155/2016/9352725]
\newblock G. A. Ngwa and M. I. Teboh-Ewungkem,
\newblock A mathematical model with quarantine states 
for the dynamics of Ebola virus disease in human populations,
\newblock {\emph{Comput. Math. Methods Med.}} {\textbf{2016}} (2016), Art. ID 9352725, 29~pp.

\bibitem{crimean} [10.1002/jmv.24312]
\newblock A. Papa, K. Tsergouli, D. \c{C}a\u{g}lay\i{k}, S. Bino, N. Como, Y. Uyar, et~al.,
\newblock Cytokines as biomarkers of Crimean-Congo hemorrhagic fever,
\newblock \emph{J. Med. Virol.} \textbf{88} (2016), no.~1, 21--27.

\bibitem{Pontryagin:1962} (MR0166037)
\newblock L. S. Pontryagin, V. G. Boltyanskii, R. V. Gamkrelidze and E. F. Mishchenko,
\newblock \emph{The mathematical theory of optimal processes},
\newblock Interscience Publishers John Wiley \& Sons, Inc.\, New York-London, 1962.

\bibitem{MR3349757} (MR3349757) [10.1155/2015/842792]
\newblock A. Rachah and D. F. M. Torres,
\newblock Mathematical modelling, simulation, 
and optimal control of the 2014 Ebola outbreak in West Africa,
\newblock \emph{Discrete Dyn. Nat. Soc.} 
\textbf{2015} (2015), Art. ID 842792, 9~pp.
\newblock {\tt arXiv:1503.07396}

\bibitem{MyID:331} [10.1002/mma.3841]
\newblock A. Rachah and D. F. M. Torres,
\newblock Predicting and controlling the Ebola infection,
\newblock \emph{Math. Method. Appl. Sci.}, in press.
\newblock {\tt arXiv:1511.06323}

\bibitem{rivers2014} [10.1371/currents.outbreaks.fd38dd85078565450b0be3fcd78f5ccf]
\newblock C. M. Rivers, E. T. Lofgren, M. Marathe, S. Eubank and B. L. Lewis,
\newblock Modeling the impact of interventions 
on an epidemic of Ebola in Sierra Leone and Liberia,
\newblock Technical report, \emph{PLOS Currents Outbreaks}, 2014.

\bibitem{rowe} [10.1086/514318]
\newblock A. K. Rowe et al., 
\newblock Clinical, virologic, and immunologic follow-up of convalescent 
Ebola hemorrhagic fever patients and their household contacts, Kikwit, Democratic Republic of the Congo,
\newblock {\emph{J. Infect. Dis.}} {\textbf{179}} (1999), Suppl.~1, S28--S35.

\bibitem{SilvaTorres:TBHIV:2015} (MR3392642) [10.3934/dcds.2015.35.4639]
\newblock C. J. Silva and D. F. M. Torres,
\newblock A TB-HIV/AIDS coinfection model and optimal control treatment,
\newblock \emph{Discrete Contin. Dyn. Syst.} \textbf{35} (2015), no.~9, 4639--4663.
\newblock {\tt arXiv:1501.03322}

\bibitem{Walsh} [10.1371/journal.pbio.0030371]
\newblock P. D. Walsh, R. Biek and L. A. Real,
\newblock Wave-like spread of Ebola Zaire,
\newblock \emph{PLoS Biology} \textbf{3} (2005), no.~11, 1946--1953.

\end{thebibliography}
\end{document}